\documentclass{article}
\usepackage{graphicx}
\def\putover#1{\mathop{\vbox{\ialign{##\crcr\noalign{\kern0pt} 
$\hfil\displaystyle{#1}\hfil$\crcr}}}\limits}
\def\ei{e(w^i)}\def\ej{e(w^j)}
\def\ek{e(w^k)}\def\el{e(w^\ell)}\def\epp{e(w^p)}\def\eijk{e(w^{ijk})}\def\eijl{e(w^{ij\ell})}\def\eijm{e(w^{ijm})}\def\eikl{e(w^{ik\ell})}\def\eikm{e(w^{ikm})}\def\eilm{e(w^{i\ell 
m})}\def\ejkl{e(w^{jk\ell})}\def\ejkm{e(w^{jkm})}\def\ejlm{e(w^{j\ell m})}
\def\eklm{e(w^{k\ell m})}
\newtheorem{tm}{Theorem}
\newtheorem{lm}{Lemma}
\newtheorem{rem}{Remark}
\newtheorem{prop}{Proposition}

\newtheorem{ex}{Example}

\newcommand   {\re}         {{\bf R}}

\newcommand   {\en}         {{\bf n}}

\newcommand   {\Yc}  {{\cal Y}}
\newcommand   {\ep}  {\varepsilon}
\newcommand    {\Uc}  {{\cal U}}
\begin{document}
\begin{center}
{\large FILTRATION LAW FOR POLYMER FLOW THROUGH POROUS MEDIA}\\
\ \\
{\bf Alain Bourgeat $\;\;\;\;\;\;\;\;\;\;$ Olivier Gipouloux}\\
Facult\'{e} des sciences et techniques, 
Universit\'{e} de St-Etienne\\
23 Rue Dr Paul Michelon, 
42 023 St-Etienne Cedex 2, France\\
\ \\
{\bf Eduard Maru\v{s}i\'{c}-Paloka}\\
Department of Mathematics, University of Zagreb,\\
 Bijeni\v{c}ka 30, 10000 Zagreb, Croatia
\end{center}
 \voffset=0in

  \setcounter{page}{0}
  \pagenumbering{arabic}

{\small{\bf Abstract.} In this paper we study the 
filtration laws for the polymeric flow in a porous 
medium. We use the quasi-Newtonian models with share 
dependent viscosity obeying the power-law and the 
Carreau's law. Using the method of homogenization 
in \cite{ad:BM} the coupled micro-macro homogenized 
law, governing the quasi-newtonian flow in a periodic 
model of a porous medium, was found. We decouple that 
law separating the micro from the macro scale. We write
the macroscopic filtration law in the form of non-linear
 Darcy's law and we prove that the obtained law is well
 posed. We give the analytical as well as the numerical
 study of our model.} 
\section{Introduction}
It is well-known that the viscosity of a polymer melt
 or a polymer solution significantly changes with the
 share rate (see e.g. \cite{BAH}).
Therefore the Newtonian 
model, characterized by a constant viscosity, is not a 
reasonable choice for modelling the polymeric flow. The 
simplest way to overcome that difficulty is to use the 
notion of quasi-Newtonian flow, i.e. to modify the 
Newton's model by using the variable viscosity that depends
 on the share-rate according to some empirical law. In this
 paper, in order to describe the polymer flow, we use the
share dependent viscosities obeying the power-law 
\begin{equation}
\eta(e(u))=\mu |e(u)|^{r-2} e(u)\;\;,\;1<r<2\;\;,\;\mu>0
\end{equation}
or the Carreau's law 
\begin{equation}
\eta(e(u)) = (\eta_0 -\eta_\infty ) (1+\lambda |e(u) |^2)^{r/2 -1} + 
\eta_\infty\;\;,\;1<r<2\;\;,\;\eta_0>\eta_\infty\geq 0\;\;,
\;\lambda>0\;,
\end{equation} 
where
\begin{eqnarray*}
&&e(u)=\frac{1}{2}(\nabla u +\nabla u^t)\;\;\mbox{the rate-of-strain 
tensor}\;\;,\\
&&|e(u)|=(\sum_{i,j=1}^n |e(u)_{ij}|^2)^{1/2} \;
-\mbox{the share rate}\;\;,
\end{eqnarray*}
To find the effective law describing the polymer flow through a 
porous medium $\Omega_\varepsilon $\footnote{where $\ep$ 
stands for the pore size} Bourgeat and Mikeli\'{c} in \cite{ad:BM} 
(see also \cite{Ho} and 
\cite{lec}) used the porous medium with periodic geometry and the 
homogenization technique called the two-scale convergence.
Starting from the microscopic problem
\begin{eqnarray}
&&-2\mbox{div}\{\eta (e(u^\ep))e(u^\ep)\} + \nabla p^\ep =
 \ep^\gamma\;f\;\;\mbox{in}\;\;\Omega_\varepsilon\\
&&\mbox{div}\,u^\ep =0 \;\;\mbox{in}\;\;\Omega_\varepsilon\\
&&u^\varepsilon =0 \;\;\mbox{on}\;\;\partial\Omega_\varepsilon\;,
\end{eqnarray}
depending on $\gamma$, the following three types of homogenized 
problems were rigorously derived:
\begin{enumerate}
\item If the Carreau's law was taken for the viscosity and
 $\gamma\neq 1$, the homogenized law is the Darcy's law
\begin{eqnarray}
&&v(x)=K(f(x)-\nabla_x\; p^0 (x))\;\;\;,\;\;x\in\Omega \label{fl,1}\\
&&\mbox{div}_x\;\;v=0\;\;\mbox{in}\;\;\Omega\;\;,
\end{eqnarray}
where $\Omega $ is the whole domain (the fluid and the solid part), 
$v$ is the filtration velocity, $p^0$ is the pressure and $K$ is the
 permeability tensor,  depending on the pore structure. In 
real-world applications $K$ is a measurable quantity. In our 
periodic model it can be computed from the cell problem
\begin{eqnarray}
&&-\mu \Delta w^j + \nabla \pi^j = e_j\;\;\mbox{in}\;\;\Yc\nonumber\\
&&\mbox{div}\,w^j =0\;\;\mbox{in}\;\;\Yc\label{locall,1}\\
&&(w^j,\pi^j)\;\;\mbox{is}\;\Yc -\mbox{periodic}\nonumber
\end{eqnarray}
by 
\begin{equation}
K_{ij} = \int_\Yc w^i_j=\mu \int_\Yc \nabla w^i \nabla w^j\;\;,
\end{equation}
where $\Yc \subset ]0,1[^n $ is the fluid part of the unit 
cell (the period) $Y=]0,1[^n $ and $S$ is the boundary of
 its solid part $A=Y\backslash \Yc$ of $Y$. The viscosity
 $\mu $ is equal to $\eta_0$, in case $\gamma <1$ and 
to $\eta_\infty $ if $\gamma>1$ and $\eta_\infty>0$. Defining
$$u^0(x,y)=\sum_{i=1}^n w^i(y)\;
\frac{\partial p^0}{\partial x_i}(x)\;\;,\;\;p^1(x,y)=
\sum_{i=1}^n \pi^i(y)\;\frac{\partial p^0}{\partial x_i}(x)\;\;,$$
the above filtration law (\ref{fl,1}) and the cell problem 
(\ref{locall,1}) can be written in the coupled form
\begin{eqnarray}
&&-\mu\;\Delta_y\;u^0+\nabla_y p^1  = f(x) - \nabla_x p^0 (x) 
\;\;\mbox{in}\;\;\Omega \times \Yc\\
&&\mbox{div}_y u^0 =0 \;\;\mbox{in}\;\;\Omega \times \Yc\\
&&(u^0,p^1)\;\;\mbox{is}\;\Yc -\mbox{periodic in }\;y\\
&&\mbox{div}_x \left(\int_\Yc u^0\, dy \right) =0\;\;\mbox{in}\;
\;\Omega\\
&&u^0 =0 \;\;\mbox{on}\;\Omega\times S\\
&&\en \cdot  \left(\int_\Yc u^0\, dy \right) =0 \;\;\mbox{on}\;\partial \Omega \;\;\;.
\end{eqnarray}
Two functions $u^0$ and $p^1$ are, in fact, first oscillating terms in the asymptotic expansion for the velocity and for the pressure, respectively.
\item If the viscosity obeys the Carreau's law with $\gamma=1$, the homogenized law is the coupled Carreau's law
\begin{eqnarray}
&&-\mbox{div}_y \{\eta_c (e_y(u^0)) e_y(u^0)\}+\nabla_y p^1  = f(x) - \nabla_x p^0 (x) \;\;\mbox{in}\;\;\Omega \times \Yc\nonumber\\
&&\mbox{div}_y u^0 =0 \;\;\mbox{in}\;\;\Omega \times \Yc\\
&&(u^0,p^1)\;\;\mbox{is}\;\Yc -\mbox{periodic in }\;y\nonumber\\
&&\mbox{div}_x \left(\int_\Yc u^0\, dy \right) =0\;\;\mbox{in}\;\;\Omega\nonumber\\
&&u^0 =0 \;\;\mbox{on}\;\Omega\times S\nonumber\\
&&\en \cdot  \left(\int_\Yc u^0\, dy \right) =0 \;\;\mbox{on}\;\partial \Omega \;\;,\nonumber\\
&&\eta_c(\xi)=(\eta_0 -\eta_\infty ) (1+\lambda |\xi|^2)^{r/2 -1} + \eta_\infty\;\;.
\end{eqnarray}
Here $p^0$ is the pressure and 
$$v=\int_{\cal Y} u^0(x,y)\;dy $$
is the filtration velocity. Functions $u^0(x,y)$  and $p^1(x,y)$
 are oscillating in $y$ and, in analogy with the linear case,
 can be seen as the first oscillating terms in the asymptotic 
expansion of the flow.\\
\item In case of power-law, 
the homogenized law is the coupled power-law
\begin{eqnarray}
&&-\mbox{div}_y \{\eta_p (e_y(u^0)) e_y(u^0)\}+\nabla_y p^1 
 = f(x) - \nabla_x p^0 (x) \;\;\mbox{in}\;\;\Omega 
\times \Yc\nonumber\\
&&\mbox{div}_y u^0 =0 \;\;\mbox{in}\;\;\Omega \times \Yc\nonumber\\
&&(u^0,p^1)\;\;\mbox{is}\;\Yc -\mbox{periodic in }\;y\label{coup,1}\\
&&\mbox{div}_x \left(\int_\Yc u^0\, dy \right) =0\;\;\mbox{in}\;
\;\Omega\nonumber\\
&&u^0 =0 \;\;\mbox{on}\;\Omega\times S\nonumber\\
&&\en \cdot  \left(\int_\Yc u^0\, dy \right) =0 \;\;\mbox{on}
\;\partial \Omega \;\;,\nonumber\\
&&\eta_p (\xi)=\mu |\xi |^{r-2}\xi\;\;.
\end{eqnarray}
\end{enumerate}
In case 1 two scales are separated and the filtration law 
is completely macroscopic. All the microscopic information
 are contained in the permeability tensor $K$, which can be measured in applications. In cases 2 and 3 the 
situation is different and two scales are still coupled.
Even if, mathematically speaking, those results solve the 
homogenization problem for such fluids, to get a physically 
relevant result one needs to decouple the problems and to 
find the filtration laws in the usual macroscopic form
\begin{equation}
v = {\cal U} (f - \nabla_x p^0 )\label{darcy,1}
\end{equation}
where 
$$v(x) = \int_\Yc u^0 (x,y)\,dy\;\;.$$
Using the idea from \cite{ad:BMP1}, \cite{ad:BMP2} and 
\cite{ad:MPM} to decouple the homogenized problems we 
define the auxiliary problem
\begin{eqnarray}
&&-\mbox{div}_y \{\eta_\alpha (e_y(w_\xi)) e_y(w_\xi)\}+
\nabla_y \pi_\xi  = \xi 
\;\;\mbox{in}\;\; \Yc\label{osnovna,1}\\
&&\mbox{div}_y w_\xi =0 \;\;\mbox{in}\;\; \Yc\\
&&(w_\xi,\pi_\xi)\;\;\mbox{is}\;\Yc -\mbox{periodic }\\
&&w_\xi =0 \;\;\mbox{on}\;S\;\;,\;\;\alpha=c,p\;\;.\label{osnovno,1}
\end{eqnarray}
With $(w_\xi,\pi_\xi)$ we define the permeability function
${\cal U} : \re^n \to \re^n $ by
\begin{equation}
{\cal U} (\xi ) = \int_\Yc w_\xi (y) \,dy\;\;,\label{jup,1}
\end{equation}
and we obtain formally the filtration law (\ref{darcy,1}). 
From the definition we can obviously conclude that the function 
${\cal U}$ is odd, i.e. that
$${\cal U} (-\xi ) = -{\cal U} (\xi )\;\;.$$
From (\ref{coup,1}) we also get
\begin{equation}
\mbox{div}\,v=0\;\;\mbox{in}\;\;\Omega\;\;,\;\;
v\cdot {\bf n}=0\;\;\mbox{on}\;\;\partial\Omega\;.\label{ellip,1}
\end{equation}
 To show that (\ref{coup,1}) and (\ref{darcy,1})-(\ref{ellip,1})
 are equivalent, we prove in  sections \ref{secar,1} and 
\ref{sepl,1} that ${\cal U}$ is monotone and coercive, 
consequently, that (\ref{darcy,1})-(\ref{ellip,1})  is well posed.\\
 We also perform the qualitative analysis of the permeability 
function ${\cal U}$. In case of coupled Carreau's law, we find 
 its Taylor's expansion.
\section{The Carreau's coupled homogenized problem}\label{secar,1}
In this section we use the notation $\eta $ for the Carreau's viscosity $\eta_c $, i.e. we omit the index $c$.
 We begin this section with statement of the main results\\
\section{Statement of the main result}
To prove that our formal separation of the coupled problem is
 meaningful we have to prove that our macroscopic problem has
 a unique solution.
\begin{tm}
Let ${\cal U}:{\bf R}^n\to{\bf R}^n $ be defined by (\ref{jup,1}).
 Then the macroscopic problem
\begin{eqnarray}
&&\mbox{div}\;{\cal U}(f-\nabla p^0)=0\;\;\;\mbox{in}\;\;\Omega\\
&& {\bf n}\cdot{\cal U}(f-\nabla p^0)=0\;\;\;\mbox{on}\;\;
\partial\Omega
\end{eqnarray}
has a unique (up to a constant) solution $p^0\in W^{1,r'} (\Omega )$,
 for any $f\in L^{r'}(\Omega)^n$.
\end{tm}\label{te,1}
In order to approximate the permeability function ${\cal U}$ 
in vicinity of $0$ we compute its Taylor's expansion. First 
two terms can be written in the form
\begin{equation}
{\cal U} (\xi ) = K\xi +\frac{1}{2} (\eta_0-\eta_\infty) 
\lambda (2-r) \sum_{m,j,k,\ell=1}^n 
H^{\ell m}_{jk} \xi_j\xi_k\xi_\ell \;e_m + O(|\xi|^5)\;,\label{474,1}
\end{equation}
where $K$ is the classical Darcy's permeability and the
coefficients $H^{\ell m}_{jk} $ are defined by (\ref{H,1}). S
uch approximate filtration law, given by the polynomial permeability
$${\cal V}(\xi)=K\xi +\frac{1}{2} (\eta_0-\eta_\infty) 
\lambda (2-r) \sum_{m,j,k,\ell=1}^n 
H^{\ell m}_{jk} \xi_j\xi_k\xi_\ell \;e_m $$
is still well posed:
\begin{tm}
Let $f\in L^4 (\Omega )^n\;\;$.
The problem 
\begin{eqnarray}  
&&\mbox{div}\;{\cal V}(f-\nabla q^0)=0\;\;\;\mbox{in}\;\;
\Omega\label{Vek,1}\\
&& {\bf n}\cdot{\cal V}(f-\nabla q^0)=0\;\;\;\mbox{on}\;
\;\partial\Omega\label{avac,1}
\end{eqnarray}
 has a  unique (up to a constant) solution 
$q^0\in W^{1,4} (\Omega )$. \label{cic,1}
\end{tm}
\subsection{Taylor's expansion of the permeability function}
We want to find the Taylor's expansion for ${\cal U}$ in 
vicinity of $0$. Since the function is odd all the 
derivatives of a pair order 
$D^\alpha {\cal U} (0)\;,\; |\alpha | =2m\;,\;m\in {\bf N}$ 
are equal to $0$. Deriving (\ref{osnovna,1}) with respect 
to $\xi_j$ we find
\begin{tm}
Let $K$ be the Darcy's permeability tensor computed from the 
local problem
\begin{eqnarray}
&&-\eta_0 \Delta w^j + \nabla \pi^j = e_j\;\;\mbox{in}\;\;
\Yc\nonumber\\
&&\mbox{div}\,w^j =0\;\;\mbox{in}\;\;\Yc\label{local,1}\\
&&(w^j,\pi^j)\;\;\mbox{is}\;\Yc -\mbox{periodic}\nonumber
\end{eqnarray}
by 
\begin{equation}
K_{ij} = \int_\Yc w^i_j=\eta_0 \int_\Yc \nabla w^i \nabla w^j\;\;.
\end{equation}
Then
\begin{equation}
\nabla_\xi {\cal U} (0) = K
\end{equation}\label{tm1,1}
\end{tm}
To prove that we proceed as in \cite{ad:BMP1}, \cite{ad:BMP2}
 or \cite{ad:MPM}. We begin by recalling the result from 
\cite{ad:sandri}:
\begin{lm}
There exist a constant $c^r_1  > 0$ such that :
\begin{equation}
c^r_1 \frac{|e(v-u)|^2_{L^r (\Yc)}}{1+|e(u)|^{2-r}_{L^r (\Yc)} 
+|e(v)|^{2-r}_{L^r (\Yc)}} 
\leq \int_{\Yc} [\eta (e(u)) e(u) -\eta (e(v)) e(v)]\;e(u-v)
\label{sandri1,1}
\end{equation}
for any $$u,v \in W=\{\phi\in W^{1,r} (\Yc)^n\;;\;
\mbox{div}\,\phi =0 \;,\;\phi \;\mbox{is}\;\Yc-\mbox{periodic}\;,
\;\phi=0\;\mbox{on}\;S\}\;.$$
\end{lm}
\begin{lm}
Let $(w_\xi, \pi_\xi)\in W\times L^2 (\Yc)$ be defined by 
(\ref{osnovna,1})-(\ref{osnovno,1}).
Then
\begin{equation}
|e(w_\xi)|_{L^r(\Yc)}\leq 2\frac{C^r}{c_1^r} 
|\Yc|^{1-\frac{1}{r}}\;|\xi |
\;\;,\label{aprior,1}
\end{equation}
for any $\xi$ such that 
$$|\xi |<\frac{c_1^r}{2C^r|\Yc|^{1-\frac{1}{r}}} \;,$$
where $C^r$ is the Poincar\'{e}-Korn's constant in $W$ 
(i.e. such that $|\phi|_{L^r(\Yc)}\leq C^r |e(\phi)|_{L^r(\Yc)} 
\;\forall \;\phi \in W $).
\end{lm}
{\bf Proof.} 
Multiplying (\ref{osnovna,1}) by $w_\xi$ and integrating over 
$\Yc$ we obtain
$$\xi\cdot \int_{\Yc} w_\xi = \int_{\Yc}\eta(e(w_\xi)) e(w_\xi)
 \geq c_1^r \frac{|e(w_\xi)|_{L^r(\Yc)}^2}{1+ 
|e(w_\xi)|_{L^r(\Yc)}^{2-r}}\;.$$
Now (\ref{sandri1,1}) gives
$$c_1^r |e(w_\xi)|_{L^r(\Yc)}\leq 
C^r|\Yc|^{1-\frac{1}{r}}(1+|e(w_\xi)|_{L^r(\Yc)}^{2-r} )\;|\xi |\;.$$
If $|e(w_\xi)|_{L^r(\Yc)}<1$ then (\ref{aprior,1}) obviously holds. 
On the other hand, supposing that $|e(w_\xi)|_{L^r(\Yc)}>1$, implies
$$c_1^r |e(w_\xi)|_{L^r(\Yc)}\leq  
2 C^r|\Yc|^{1-\frac{1}{r}}|e(w_\xi)|_{L^r(\Yc)}^{2-r} \;|\xi |\;$$
and
$$|e(w_\xi)|_{L^r(\Yc)}^{r-1}\leq  
2 \frac{C^r}{c_1^r}|\Yc|^{1-\frac{1}{r}} \;|\xi |\;$$
which contradicts the assumption that 
$$|\xi |<\frac{c_1^r}{2C^r|\Yc|^{1-\frac{1}{r}}} \;.\;\clubsuit$$
{\bf Proof of theorem 3.}
To prove the claim we only need to show that, as $|h|\to 0$
\begin{equation}
|w_h - W(h) |_W = o(|h|)\;\;,\;\;W(h)=\sum_{i=1}^n h_i\, w^i
\end{equation}
because
\begin{equation}
|{\cal U} (h) - K\,h |=|\int_{\Yc} [w_h -W(h)] |\leq 
C |w_h - W(h)|_W \;.
\end{equation}
In fact we shall prove that
\begin{equation}
|{\cal U} (h) - K\,h |=O(|h|^3)
\end{equation}
which implies, not only that $\nabla {\cal U} (0) =K $, but 
also that $\nabla^2 {\cal U} (0) =0$ (which we knew before 
because ${\cal U}$ is an odd function). Using the standard 
regularity result (see e.g. \cite{ad:KF}) for the problem 
(\ref{local,1}) we see that
\begin{equation}
|W(h)|_{H^3 (\Yc )} \leq \sum_{i=1}^n |h|\;|w^i|_{H^3 (\Yc )}\leq 
C\,|h|\;\;\;.\label{h2,1}
\end{equation}
and therefore
\begin{equation}
|W(h)|_{C^1 (\overline{\Yc} )} \leq C\,|h|\;\;\;.\label{c1,1}
\end{equation}
Functions $W(h)$ and $\Pi (h) = \sum_{i=1}^n h_i\, \pi^i $ 
satisfy the system
\begin{eqnarray}
&&-\mbox{div}\{\eta \{e[W(h)]\} e[W(h)]\} +
 \nabla \Pi (h) = h + R_h\label{32,1}\\
&&\mbox{div}\,W(h) =0\\
&&(W(h) ,\Pi (h) )\;\;\mbox{is}\;\Yc-\mbox{periodic}\;,\;W(h) =0 
\;\;\mbox{on}\;S\;,
\end{eqnarray}
where
$$R_h =\mbox{div}\{(\eta_0 -\eta_\infty) 
[(1+\lambda |e[W(h)]|^2)^{r/2\;-1} -1] e[W(h)] \}\;\;.$$
An easy computation leads to
$$(1+\lambda t^2)^{r/2 -1} -1 = 
(r/2 -1) t^2 \int_0^\lambda (1+\tau t^2)^{r/2 -2} 
d\tau \leq C t^2 \;\;.$$
For $1<r<2$ we get using (\ref{c1,1})
\begin{equation}
\int_{\Yc} R_h\,\phi \leq C |\phi |_{W^{1,r}(\Yc)}|\;
|e[W(h)]|^3_{L^{3r'}(\Yc)}\leq C |h|^3 
|\phi |_{W^{1,r}(\Yc)}\;.\label{rh2,1}
\end{equation}
Subtracting (\ref{osnovna,1}) from (\ref{32,1}), multiplying by 
$W(h) -w_h$ and integrating over $\Yc$ we obtain
\begin{eqnarray*}
&&J=\int_{\Yc} \{ \eta\{e[W(h)]\} e[W(h)] - \eta\{e(w_h)\} e(w_h)\} 
\;e(W(h)-w_h) =\\
&&=\int_{\Yc} R_h (W(h)-w_h) \leq C |h|^3 |W(h)-w_h |_W\;\;.
\end{eqnarray*}
An application of lemma 1 gives 
$$J\geq c^r_1 \frac{|W(h)-w_h |^2_W}{1 + C|h|} \;\;.\; \clubsuit $$
\ \\
As we have seen the second derivative in $\xi=0$ is $0$.  Third 
can be computed from the auxiliary problem
\begin{eqnarray}
&&-\eta_0 \Delta w^{\ell jk} + \nabla \pi^{\ell jk} = \\
&&=(\eta_0-\eta_\infty ) (r-2) \lambda 
\mbox{div}\{ e(w^k)\cdot e(w^j)\; e(w^\ell) +\\
&&+ e(w^k)\cdot e(w^\ell)\; e(w^j)+e(w^\ell)\cdot e(w^j)\; 
e(w^k)\}\;\mbox{in}\;\Yc\\
&&\mbox{div}\,w^{ijk} =0\;\;\mbox{in}\;\;\Yc\\
&&(w^{ijk},\pi^{ijk})\;\;\mbox{is}\;\Yc -\mbox{periodic}
\end{eqnarray}
by
$$\frac{\partial^3 {\cal U}(0)}{\partial \xi_\ell \partial \xi_j 
\partial\xi_k} = \int_\Yc w^{\ell jk}\;\;.$$
Since
\begin{eqnarray}
&&\frac{\partial^3 {\cal U}(0)_m}{\partial \xi_\ell \partial \xi_j 
\partial\xi_k} =\eta_0 \int_\Yc \nabla w^{\ell j k} \nabla w^m =\\
&&=(\eta_0-\eta_\infty ) (2-r)\lambda \int_\Yc \{ e(w^k)\cdot e(w^j)
\; e(w^\ell)\cdot e(w^m) +\\
&&+ e(w^k)\cdot e(w^\ell)\; e(w^j)\cdot e(w^m)+
e(w^\ell)\cdot e(w^j)\; e(w^k)\cdot e(w^m)\}=\\
&&=(\eta_0-\eta_\infty ) (2-r)\lambda (H^{\ell m}_{jk} + 
H^{jm}_{\ell k} + H^{km}_{\ell j} )\;\;,
\end{eqnarray}
where
\begin{equation}
H^{\ell m}_{jk} = H^{m\ell}_{jk}=H^{\ell m}_{kj} =H^{jk}_{\ell m} 
=\int_\Yc  e(w^k)\cdot e(w^j)\; e(w^\ell)\cdot e(w^m)\;\;.
\label{H,1}
\end{equation}
This gives the expansion for ${\cal U}$ in the form
\begin{equation}
{\cal U} (\xi ) = K\xi +\frac{1}{2} (\eta_0-\eta_\infty) 
\lambda (2-r) \sum_{m,j,k,\ell=1}^n 
H^{\ell m}_{jk} \xi_j\xi_k\xi_\ell \;e_m + O(|\xi|^5)\;.\label{47,1}
\end{equation}
\begin{prop}
Let $H^{\ell m}_{jk} $ be defined by (\ref{H,1}). Then 
(\ref{47,1}) holds.
\end{prop}
{\bf Proof.}
We use the same idea as in the proof of theorem \ref{tm1,1}. We
 define
\begin{eqnarray*}
&&V(h)=\sum_{i+1}^n w^i h_i +\frac{1}{2}
\sum_{i,j,k=1}^n w^{ijk} h_i h_j h_k\\
&&Q(h)=\sum_{i+1}^n \pi^i h_i +\frac{1}{2}\sum_{i,j,k=1}^n 
\pi^{ijk} h_i h_j h_k\;\;.
\end{eqnarray*}
Now
\begin{eqnarray*}
&&-\mbox{div}\{\eta_r(e[V(h))]e(V(h))\}+\nabla Q(h)=\\
&&=h-(\eta_0-\eta_\infty )\mbox{div}\{[1+
\lambda |e(V(h))|^2)^{r/2-1} -1] e[V(h)]\}+\\
&&+\frac{\lambda}{2} (r-2)(\eta_0-\eta_\infty )
\mbox{div}\{\sum_{i,j,m=1}^n e(w^i) \cdot e(w^j)\;
e(w^m)\;h_i h_j h_m \}\;\;.
\end{eqnarray*}
It only remains to estimate
\begin{eqnarray*}
&&J=(\eta_0-\eta_\infty )
\mbox{div}\{[1+\lambda |e(V(h))|^2)^{r/2-1} -1] 
e[V(h)]\}+\\
&&+\frac{\lambda}{2} (r-2)(\eta_0-\eta_\infty )
\mbox{div}\{\sum_{i,j,m=1}^n e(w^i) \cdot e(w^j)\;
e(w^m)\;h_i h_j h_m \}\;\;.
\end{eqnarray*}
A simple Taylor's formula gives
$$(1+\lambda |\xi |^2 )^{r/2 -1} - 1 =
\frac{1}{2}\lambda |\xi |^2 + O(\lambda^2 |\xi|^4)\;\;.$$
After recalling that 
$e(w^i)\in C(\overline{\cal Y})^{n\times n}\;,
\;e(w^{ijk})\in C(\overline{\cal Y})^{n\times n} $, we
 get for any $y\in {\cal Y}$
\begin{eqnarray*}
&&(\eta_0-\eta_\infty )\mbox{div}\{[1+
\lambda |e(V(h))|^2)^{r/2-1} -1] e[V(h)]\}=\\
&&=\frac{\lambda}{2} (r-2)(\eta_0-\eta_\infty )
\mbox{div}\{\sum_{i,j,m=1}^n e(w^i) \cdot e(w^j)\;
e(w^m)\;h_i h_j h_m \}+O(|h|^4)\;\;.
\end{eqnarray*}
 Arguing as in the proof of theorem \ref{tm1,1}, we get 
the claim.$\;\clubsuit$\\
\begin{rem}
It can be seen from the definition that 
$H^{ii}_{jj}>0\;,\;H^{\ell k}_{\ell k}>0$. Furthermore, 
for
$$Z(h)=\sum_{i,j,k,\ell=1}^n 
H^{i\ell}_{jk} h_j h_k h_i \;{\bf e}_\ell=
\int_{\cal Y} |\sum_{k=1}^n e(w^k) h_k |^2 
\sum_{i=1}^n e(w^i) h_i \sum_{\ell=1}^n e(w^\ell)\;{\bf e}_\ell $$
we have
$$Z(h)\cdot h=\int_{\cal Y}|\sum_{k=1}^n e(w^k) h_k|^2 > 0\;\;.$$
\end{rem}
To compute further order terms we proceed in the same way and we 
get another auxiliary problem
\begin{eqnarray}
&&-\eta_0 \Delta w^{ijk\ell p} + \nabla \pi^{ijk\ell p} = \\
&&= (\eta_0 - \eta_\infty)\lambda (r-2)\mbox{div}\{\\
&&\ei\cdot\ej\eklm +\ek\cdot\ei\ejlm +\ei\cdot\el\ejkm +\\ 
&&\ei\cdot\epp\ejkl +\ej\cdot\ek\eilm +\ej\cdot\el\eikm +\\ 
&&\ej\cdot\epp\eikl +\ek\cdot\el\eijm +\ek\cdot\epp\eijl +\\ 
&&\el\cdot\epp\eijk +  \\ 
&&(\ejkl\cdot\ei+\eikl\cdot\ej+\eijl\cdot\ek+\eijk\cdot\el )\epp +\\
&&(\ejkm\cdot\ei+\eikm\cdot\ej+\eijm\cdot\ek+\eijk\cdot\epp )\el +\\
&&(\ejlm\cdot\ei+\eilm\cdot\ej+\eijm\cdot\el+\eijl\cdot\epp )\ek +\\
&&(\eklm\cdot\ei+\eilm\cdot\ek+\eikm\cdot\el+\eikl\cdot\epp )\ej +\\
&&(\eklm\cdot\ej+\ejlm\cdot\ek+\ejkm\cdot\el+\ejkl\cdot\epp )\ei +\\
&&\lambda(r-4)\{ \\
&&\{\el\cdot\ek\ej\cdot\ei+\el\cdot\ej\ek\cdot\ei+\el\cdot\ei\ek\cdot
\ej\}\epp +\\
&&\{\epp\cdot\ek\ej\cdot\ei+\epp\cdot\ej\ek\cdot\ei+\epp\cdot\ei\ek
\cdot\ej\}\el +\\
&&\{\epp\cdot\el\ej\cdot\ei+\epp\cdot\ej\el\cdot\ei+
\epp\cdot\ei\el\cdot\ej\}\ek +\\
&&\{\epp\cdot\el\ek\cdot\ei+\epp\cdot\ek\el\cdot\ei+
\epp\cdot\ei\el\cdot\ek\}\ej +\\
&&\{\epp\cdot\el\ek\cdot\ej+\epp\cdot\ek\el\cdot\ej+
\epp\cdot\ej\el\cdot\ek\}\ei \}\}\\
&&\mbox{div}\,w^{ijk\ell p} =0\;\;\mbox{in}\;\;\Yc\\
&&(w^{ijk\ell p},\pi^{ijk\ell p})\;\;\mbox{is}\;\Yc -\mbox{periodic .}
\end{eqnarray}
Now
$$\frac{\partial^5 {\cal U}(0)}{\partial \xi_i \partial \xi_j 
\partial \xi_k 
\partial \xi_\ell \partial \xi_p}=\int_{\cal Y} w^{ijk\ell p} $$
and we get an expansion in the form
\begin{eqnarray*}
&&{\cal U}(\xi)=K\xi+\frac{1}{2}(\eta_0 -\eta_\infty ) \lambda (2-r) 
\sum_{i,j,k,\ell,m=1}^n H^{\ell m}_{jk} \xi_j \xi_k \xi_\ell\;
{\bf e}_m + \\
&&+\frac{1}{2}(\eta_0 -\eta_\infty )\lambda (2-r) 
\sum_{j,k,\ell ,m,p,q=1}^n H^{mpq}_{jk\ell} \xi_j \xi_k 
\xi_\ell \xi_m \xi_p\;{\bf e}_q + O(|\xi|^7)
\end{eqnarray*}
where
$$ H^{mpq}_{jk\ell}=\int_{\cal Y}\{e(w^j)\cdot e(w^k)\;
e(w^{\ell mp})\cdot e(w^q )+ 
\lambda \frac{(r-4)}{2} e(w^j)\cdot e(w^k)\;e(w^\ell )\cdot e(w^m )
\;e(w^p )\cdot e(w^q)\}\;.$$
This computation can be proceeded.\\
\ \\
\subsection{Proofs of existence theorems}To prove the theorem 
\ref{te,1} we first prove the following lemma that gives 
coercivity and boundness of the differential operator 
$$A(\phi)=-\mbox{div}\,{\cal U}(f-\nabla \phi)\;\;.$$
\begin{lm}
There exist constants  $M,C,c>0$ such that
\begin{eqnarray}
&&{\cal U}(\xi)\cdot\xi \geq c|\xi |^{r'}\;\;\;\mbox{for any}
\;\;|\xi|>M\label{lower,1}\\
&&|{\cal U} (\xi )|\leq C |\xi |^{r'-1}\;\;\mbox{for any}\;\;\xi\in {\bf 
R}^n\;.\label{upper,1}
\end{eqnarray}
\end{lm}\label{lele,1}
{\bf Proof.}
The upper bound is easy to prove, since we have, due to (\ref{sandri1,1})
$$c_1^r \frac{|e(w_\xi )|_{L^r ({\cal Y})}}{|e(w_\xi )|_{L^r ({\cal Y})}^{2-r} 
+1}\leq C|\xi |\;\;.$$
To prove the lower bound (\ref{lower,1}) we define $v_\xi=|\xi|^{1-r'}\;w_\xi 
$ and $q_\xi=|\xi|^{-1}\;\xi$. Those functions satisfy the system
\begin{eqnarray}
&&-\mbox{div}_y \{\eta_\xi (e_y(v_\xi)) e_y(v_\xi)\}+\nabla_y q_\xi  = 
\frac{\xi}{|\xi |}\;\;\mbox{in}\;\; \Yc\label{osno,1}\\
&&\mbox{div}_y v_\xi =0 \;\;\mbox{in}\;\; \Yc\\
&&(v_\xi,q_\xi)\;\;\mbox{is}\;\Yc -\mbox{periodic }\\&&v_\xi =0 
\;\;\mbox{on}\;S\;\;\;,\label{osnova,1}
\end{eqnarray}
where 
$$\eta_\xi (\tau )= (\eta_0 -\eta_\infty ) (|\xi |^{2(1-r')} +\lambda |\tau 
|^2 )^{r/2\;-1} +\eta_\infty |\xi |^{-r'}\;\;.$$
Applying again (\ref{sandri1,1}) we obtain
$$\frac{|e(v_\xi)|_{L^r ({\cal Y})}}{|\xi|^{2(1-r')(2-r)} +|e(v_\xi 
)|_{L^r({\cal Y })}}\leq C $$
so that we can extract a subsequence $\{\xi_n \} $ such that 
$|\xi_n| \to \infty\;,\;|\xi_n|^{-1}\;\xi_n\to \xi_\infty\;,
\;|\xi_\infty|=1 $ and
$$v_{\xi_n} \rightharpoonup v_\infty \;\;\;\;\mbox{weakly in}
\;\;\;W\;\;\mbox{as}\;\;n\to \infty\;\;.$$
Using the Minty's lemma we get that $v_\infty $ satisfies
 the system\begin{eqnarray*}
&&-(\eta_0 -\eta_\infty )\lambda^{r/2\;-1}\,\mbox{div}_y 
\{|e_y(v_\infty))|^{r-2} e_y(v_\infty)\}+\nabla_y q_\infty  = 
\xi_\infty\;\;\mbox{in}\;\; \Yc\\
&&\mbox{div}_y v_\infty =0 \;\;\mbox{in}\;\; \Yc\\
&&(v_\infty,q_\infty)\;\;\mbox{is}\;\Yc -\mbox{periodic }\\
&&v_\infty =0 \;\;\mbox{on}\;S\;\;\;,
\end{eqnarray*}
For any $\xi_\infty $, such that $|\xi_\infty | =1$, we can define 
a continuous function 
$${\cal I} (\xi_\infty)=\int_{\cal Y}v_\infty \neq 0\;\;.$$
Furthermore
$${\cal I}(\xi_\infty)\cdot \xi_\infty = 
\int_{\cal Y} (\eta_0 -\eta_\infty )\;
\lambda^{r/2\;-1}\;|e(v_\infty)|^r\;\;.$$
and there exists
$$c_\infty=\inf_{|\tau |=1} {\cal I}(\tau)\cdot \tau >0\;\;.$$
Two functionals
\begin{eqnarray*}
&&\Phi_\xi (v)=\int_{\cal Y} \eta_\xi (e(v))\;|e(v)|^2 \\
&&\Phi_\infty (v)=\int_{\cal Y} (\eta_0 -\eta_\infty )\;
\lambda^{r/2\;-1}\;|e(v)|^r
\end{eqnarray*}
are obviously convex and, due to the Lebesgues dominated 
convergence theorem,
$$\Phi_\xi (v_\xi )- \Phi_\infty (v_\xi ) \to 0\;\;\;\;\mbox{as}
\;\;\;|\xi |\to \infty\;\;.$$
Since $\Phi_\infty $ is weak lower semicontinuous in $W$, we obtain
$$\lim_{|\xi |\to \infty } \inf \Phi_\xi (v_\xi )= 
\lim_{|\xi |\to \infty }\inf  \{ \Phi_\infty (v_\xi )+
[\Phi_\xi (v_\xi )- \Phi_\infty (v_\xi )]\}\geq 
\Phi_\infty (v_\infty )\geq c_\infty >0\;\;.$$
But $$\Phi_\xi (v_\xi )= |\xi |^{-r'}\;{\cal U} (\xi )\cdot \xi $$
which proves the claim. $\;\clubsuit$\\
{\bf Proof of theorem 3.} 
Lemma \ref{lele,1} implies coercivity and boundness of the operator
 $A :W^{1,r'}(\Omega )/{\bf R} \to (W^{1,r'} (\Omega ))'$ .
 Its strict monotonicity follows from inequality
\begin{eqnarray*}
&[{\cal U}(\xi ) - {\cal U} (\tau )]\cdot (\xi -\tau)&
=\int_{\cal Y}\{\eta_c (e(w_\xi ))e(w_\xi)- \eta_c (w_\tau ))]
\cdot (w_\xi -w_\tau )\\
&&\geq c_1^r\frac{|e(w_\xi -w_\tau )|^2_{L^r ({\cal Y})}}
{|e(w_\xi )|_{L^r ({\cal Y})}^{2-r} + 
|e(w_\tau )|_{L^r ({\cal Y})}^{2-r}+1}\;\;\;.
\end{eqnarray*}
Now the classical Browder-Minty's theorem proves the existence
 of the solution, i.e. the surjecivity of our operator. As it 
is strictly monotone, it is also injective, i.e. the solution 
is unique. $\;\clubsuit $\\
\ \\
If we approximate ${\cal U}$ by taking only first two terms in its
 Taylor's series, i.e., if we consider
$${\cal V}(\xi)= K\xi +\frac{1}{2} (\eta_0-\eta_\infty) \lambda 
(2-r) \sum_{m,j,k,\ell=1}^n H^{\ell m}_{jk} \xi_j\xi_k\xi_\ell \;
e_m $$
instead of ${\cal U}$, we get the problem
(\ref{Vek,1})-(\ref{avac,1}).\\
For such approximation of the filtration law we have theorem 
\ref{cic,1}, from the beginning of this section. We are now able
 to prove it.\\
{\bf Proof of theorem \ref{cic,1}.}
We define the formal differential operator
$$B(\varphi )=\mbox{div}\,{\cal V} (f-\nabla \varphi )\;\;.$$
Since ${\cal V} $ is a polynomial of order $3$, it satisfies the estimate
$$|{\cal V} (\xi)|\leq C\;(|\xi |+|\xi |^3)\;\;.$$
Consequently $B:W^{1,4}(\Omega )/{\bf R} \to (W^{1,4}(\Omega )/{\bf R} )' $ is 
a bounded operator. Moreover
$$|B(\varphi )|_{(W^{1,4}(\Omega )/{\bf R})'} \leq C \;|f-\nabla \varphi 
|^3_{L^4 (\Omega )}\;\;.$$ 
Next we prove the coercivity. Let 
$$E(\xi )=\sum_{i=1}^n e(w^i)\;\xi_i\;\;.$$
Then, for any $\xi \in {\bf R}^n$,
$${\cal V}(\xi )\cdot \xi = \int_{\cal Y}(\;\eta_0\;
|E(\xi )|^2 +(\eta_0 -\eta_\infty )\;\lambda (1-\;r/2)|E(\xi )|^4\; )\;\;.$$
The functions
$${\cal G} (\xi)=\int_{\cal Y} |E(\xi )|^2\;\;\;\;,\;\;\;\;{\cal 
H}(\xi)=\int_{\cal Y}|E(\xi )|^4 $$
are continuous and strictly positive for $\xi\neq 0$. Therefore there exist
$$ \kappa_1=\inf_{|\xi|=1} {\cal G} (\xi 
)\;\;\;\;,\;\;\;\;\kappa_2=\inf_{|\xi|=1} {\cal H} (\xi )\;\;.$$
But than
$${\cal V}(\xi )\cdot\xi \geq \;\kappa_1\;|\xi |^2 +(\eta_0 
-\eta_\infty)\;\lambda\;(1-\;r/2)\;\kappa\;|\xi |^4 \;\;,$$
implying the coercivity of $B$ on $W^{1,4} (\Omega )/{\bf R}$. It remains to 
prove the monotonicity of $B$. Then the result will follow from the classical 
Browder-Minty's theorem.\\For any $\xi , \tau \in {\bf R}^n $ the derivative 
$D {\cal V} (\xi ) $ 
satisfies
\begin{eqnarray*}
&&D {\cal V} (\xi )\;\tau \,\cdot\, \tau =K\tau\cdot\tau+\frac{1}{2}\lambda 
(2-r) (\eta_0 -\eta_\infty) \int_{\cal Y} [(E(\xi )\cdot E(\tau ))^2 + \\ 
&& + E(\xi )^2\, E(\tau )^2 ]>K\tau\cdot\tau \geq k_0 |\tau |^2\;\;.
\end{eqnarray*}
The above inequality proves that the equation (\ref{Vek,1}) is uniformly 
elliptic , i.e. that the operator $B$ is monotone. $\;\clubsuit$\\
\ \\
The approximation ${\cal V}$ of the permeability function is valid either for 
small $|\xi |$ or for small $\lambda$. For large $\lambda $ and $|\xi |$, the 
Carreau's auxiliary problem (\ref{osnovna,1})-(\ref{osnovno,1}), with 
$\eta_c$, is close to the power-law auxiliary problem 
(\ref{osnovna,1})-(\ref{osnovno,1}), with $\eta_p$ and 
$\mu=(\eta_0-\eta_\infty)\lambda^{r/2-1}$. We study the power-law auxiliary 
problem in the next section.\\
We refer to   section \ref{car} for numerical computation of the permeability 
function ${\cal U}$ and its approximation ${\cal V}$ and their comparison  for 
sufficiently small $|\xi |$.\\
 \section{The power-law coupled homogenized problem}\label{sepl,1}
 In this section we also start with statement of its main results and we adopt 
the notation $\eta $ for the power-law viscosity $\eta_p$.
\section{Statement of the main result}
From (\ref{osnovna,1})-(\ref{jup,1}) we get by substitution that the 
permeability function is homogeneous in the sense that
\begin{equation}
{\cal U} (\lambda\,\xi ) = |\lambda|^{r'-2}\,\lambda\,{\cal U} (\xi 
)\;\;,\;\;\lambda\in {\bf R}\;\;.\label{hoe,1}
\end{equation}implying 
\begin{equation}
m |\xi |^{r'-1}\leq |\Uc (\xi )| \leq M |\xi |^{r'-1} 
\label{koe,1}\end{equation}with
$$m=\inf_{|\xi|=1} |\Uc (\xi )|\;\;,\;\;M=\sup_{|\xi|=1} |\Uc (\xi )|\;\;.$$
Now, as in the previous chapter, (\ref{koe,1}) and (\ref{mon,1}) imply that:
\begin{tm} For any $f\in L^{r'} (\Omega )^n$ the macroscopic problem
\begin{eqnarray}
&&\mbox{div}\,{\cal U}(f-\nabla p^0 )=0\;\;\;\mbox{in}\;\;\Omega\\
&&{\bf n}\cdot{\cal U}(f-\nabla p^0) =0\;\;\;\mbox{on}\;\;\partial\Omega
\end{eqnarray}
has a unique (up to a constant) solution $p^0 \in W^{1,r'} (\Omega )$.
\end{tm} The standard engineering model (see e.g. \cite{BSL}, \cite{WPW}) 
suggests that the permeability function ${\cal U} $ can be written in the form
\begin{equation}
{\cal U}(\xi)_i=|(K\,\xi)_i|^{r'-2}\;(K\,\xi)_i\;\;,\;\;i=1,\ldots 
,n\;\;,\label{conjj,1}
\end{equation}
where $K$ is a permeability tensor.\\
To check that conjecture we define the conjugated function $G$ 
by$$G(\xi)_i=|{\cal U }(\xi )_i|^{r-2} \;{\cal U}(\xi)_i\;\;,$$
such that
$${\cal U}(\xi)_i=|G(\xi)_i|^{r'-2}\,G(\xi )_i\;\;.$$
Now, if the conjecture (\ref{conjj,1}) is true, the function $G$ is linear. 
Homogeneity (\ref{hoe,1}) of ${\cal U}$ implies that $G$ is homogeneous with 
order $1$, i.e. that
$$G(\lambda\,\xi)=\lambda\,G(\xi )\;\;,\;\lambda\in {\bf R}\;\;.$$
{\em Function homogeneous with order $1$ is linear if and only if it is 
differentiable for} $\xi=0$.\\
Function $G$ is differentiable for any $\xi\neq 0$, but, in general, not for 
$\xi=0$. That means that {\bf the conjecture (\ref{conjj,1}) is, in general, 
false.} In case of unidirectional flow it is differentiable and, therefore, 
linear. In case of a porous medium consisting of a net of narrow channels, $G$ 
allows an asymptotic expansion in powers of the thickness of the channels, 
with the first term being a linear function.
Numerical examples show that $G$ is linear (or differentiable for $\xi=0$) in 
case of periodic porous medium with the solid part being an array of 
rectangular obstacles (i.e. with the fluid part being a net of perpendicular 
channels). On the other hand, in case of spherical or elliptical obstacle, the 
function $G$ is not linear and conjecture (\ref{conjj,1}) does not hold.
\subsection{Study of the permeability function}
 We next recall the result from \cite{ad:sandri}, that will be used in the 
sequel:
\begin{lm}
There exist a constant $c^r_1  > 0$ such that :
\begin{equation}
c^r_1 \frac{|e(v-u)|^2_{L^r (\Yc)}}{|e(u)|^{2-r}_{L^r (\Yc)} 
+|e(v)|^{2-r}_{L^r (\Yc)}} \leq \int_{\Yc} [\eta (e(u)) e(u) -\eta (e(v)) 
e(v)]\;e(u-v)\label{sandrii1,1}
\end{equation}
for any $$u,v \in W=\{\phi\in W^{1,r} (\Yc)^n\;;\;\mbox{div}\,\phi =0 
\;,\;\phi \;\mbox{is}\;\Yc-\mbox{periodic}\;,\;\phi=0\;\mbox{on}\;S\}\;.$$
\end{lm}
It implies that, for $\xi ,\tau \in {\bf R}^n $, $\xi\neq \tau 
$,\begin{equation}
[{\cal U}(\xi)-{\cal U}(\tau )]\cdot (\xi -\tau )\geq c^r_1 
\frac{|e(w_\xi-w_\tau)|^2_{L^r (\Yc)}}{|e(w_\xi)|^{2-r}_{L^r (\Yc)} 
+|e(w_\tau)|^{2-r}_{L^r(\Yc)}} > 0\;\;. \label{mon,1}
\end{equation}
Before we continue we give two simple examples mentioned in the introduction 
of this section.\\
\begin{ex}
As it was noticed in \cite{ad:BM}, \cite{lec}, \cite{Ho}, if the flow is 
unidirectional, e.g. in the direction ${\bf e}_1$, then ${\cal U}(\xi)={\cal 
J}(\xi_1)\;{\bf e}_1 $, with
$${\cal J}(\xi_1)=|\xi_1|^{r' -2}\;\xi_1\;\;.$$
This result corresponds to the usual models that can be found in the 
engineering literature (see e.g. \cite{WPW}, \cite{CM}).\end{ex}
\begin{ex} In case of a porous medium consisting of a network of thin pipes we 
can find the asymptotic form of our function ${\cal U}$ and, again, it 
corresponds to the usual engineering model from \cite{CM} and \cite{WPW}.\\
Suppose that $n=2$ and that ${\cal Y}={\cal Y}_\delta $ consists of two narrow 
perpendicular channels
$$P_\delta^1=]0,1[\;\times\;]\,-\delta/2 
,\delta/2\,[\;\;\;,\;\;\;P_\delta^2=\;]\,-\delta/2 
,\delta/2\,[\;\times\;]0,1[\;\;,$$
with thickness $\delta$.
We consider our auxiliary problem in such unit cell ${\cal 
Y}_\delta=P^1_\delta\cup P^2_\delta $
\begin{eqnarray}
&&-\mbox{div}_y \{\eta_p (e_y(w_\xi^\delta)) e_y(w_\xi^\delta)\}+\nabla_y 
\pi_\xi^\delta  = \xi \;\;\mbox{in}\;\; \Yc_\delta\label{osa,1}\\
&&\mbox{div}_y w_\xi^\delta =0 \;\;\mbox{in}\;\; \Yc_\delta\\
&&(w_\xi,\pi_\xi)\;\;\mbox{is}\;\Yc_\delta -\mbox{periodic }\\
&&w_\xi =0 \;\;\mbox{on}\;S^\delta\;\;,\label{os,1}
\end{eqnarray}
where $S^\delta=\partial (P^1_\delta \cup P^2_\delta )$. We study the limit as 
$\delta\to 0$. The similar problem was considered in \cite{ES} for the 
Newtonian fluid and an explicit formula for the permeability was found. We are 
going to do the same thing with our nonlinear permeability function ${\cal 
U}^\delta (\xi)=\int_{\cal Y} w_\xi^\delta $. As in \cite{ES}, we seek an 
approximation of $(w^\delta_\xi , \pi^\delta_\xi )$ in the form
$$w(\delta)=\delta^{r'}\left\{\begin{array}{l}\theta(y_2/\delta 
)\;|\xi_1|^{r'-2}\xi_1\;{\bf e}_1\;\;\mbox{in}\;\;P^\delta_1\\
\theta(y_1/\delta )\;|\xi_2|^{r'-2}\xi_2\;{\bf e}_2\;\;\mbox{in}\;\;P^\delta_2
\end{array}\right.\;\;,\;\;
\pi(\delta)=\left\{\begin{array}{l}
y_2\;\xi_2\;\;\mbox{in}\;\;P^\delta_1\\
y_1\;\xi_1\;\;\mbox{in}\;\;P^\delta_2
\end{array}\right.$$
where
$$\theta(z)=\frac{2}{r'\;\mu}\left(\frac{\sqrt{2}}{\mu}\right)^{r'-2}\;[\,(1/2)^{r'} 
-|z|^{r'} \,]\;\;.$$
In each $P^i_\delta $ we have the Poincar\'{e}-Korn's inequality
$$|w^\delta_\xi -w(\delta) |_{L^r(P^i_\delta )}\leq C\delta | e(w^\delta_\xi 
-w(\delta))|_{L^r (P^i_\delta )} $$
and the trace inequality
$$|w^\delta_\xi -w(\delta)|_{L^r (\Sigma^i_\delta )}\leq 
C\;\delta^{1/r'}\;|e(w^\delta_\xi -w(\delta)) |_{L^r (P^i_\delta )}\;\;,$$
where
$$\Sigma^1_\delta=\partial P^1_\delta \cap 
P^2_\delta\;\;\;,\;\;\;\Sigma^2_\delta=\partial P^2_\delta \cap 
P^1_\delta\;\;.$$
By a direct computation we obtain
\begin{eqnarray*}&&\mu\int_{P^i_\delta} \{|e(w^\delta_\xi )|^{r-2} 
e(w^\delta_\xi )-|e(w(\delta ))|^{r-2} e(w(\delta ))\} e(w^\delta_\xi 
-w(\delta ))=\\
&&=\int_{\Sigma^i_\delta } \{\pi (\delta ){\bf n} -\mu |e(w(\delta ))|^{r-2} 
e(w(\delta )){\bf n} \}\;(w^\delta_\xi -w(\delta))\;\;. \end{eqnarray*}
We use (\ref{sandrii1,1}) to estimate from below the left-hand side and the 
above inequalities to estimate from above the right-hand side. It leads to the 
estimate
$$|\delta^{-(1+r')}\;{\cal U}_\delta (\xi )-\langle \theta \rangle\;{\cal J} 
(\xi )\;|\leq C\;\delta^{1/r}\;\;,$$
where
$$\langle\theta\rangle=\int_{-1/2}^{1/2} \theta(s)ds 
=\frac{2}{(r'+1)\mu}(\sqrt{2}/\mu )^{r' -2}\;\;,$$
and
$${\cal J}(\xi)=(|\xi_1|^{r'-2}\xi_1 ,|\xi_2 |^{r'-2} \xi_2 )\;\;.$$
Therefore
$${\cal U}_\delta (\xi)_i=\delta^{r' +1} \frac{2}{(r'+1)\,\mu} (\sqrt{2}/\mu 
)^{r'-2}\;|\xi_i|^{r'-2}\;\xi_i + O(\delta^{r'+1+1/r})\;\;,\;\;i=1,2 $$
i.e.
$${\cal U}_\delta (\xi)\approx C_\delta J(\xi )\;\;.$$\label{e,1}
\end{ex}
We are now tempted to believe that the permeability function can be written in 
the form (\ref{conjj,1}).\\
First of all, we can not hope to have $K$ that depends only on the geometry of 
the porous medium, as in the linear case, because in two above examples it 
depends on the flow index $r$.\\
We know that this conjecture is true in case of unidirectional flow and, up to 
some order of approximation, for a network of narrow channels. \\
\ \\
To verify whether (\ref{conjj,1}) is true or not, we define the conjugated 
permeability function
$$G(\xi)_i=|\Uc (\xi)_i|^{r-2} \Uc (\xi )_i\;\;,\;\;i=1,\ldots, n\;\;.$$
As $r-2<0$ and $\Uc (0)=0$ , $G(0)$ is not correctly defined by the above 
formula. However, since $\Uc$ satisfies (\ref{hoe,1}), and (\ref{koe,1})we have
$$m |\xi |\leq |G (\xi )| \leq M |\xi |\;\;.$$
Therefore we define $G$ in $\xi=0$ by $G(0)=0$ and we get the function $G \in 
C({\bf R})^n\; $.
Function $G$ is obviously homogeneous
$$G(\lambda \,\xi)=\lambda \,G(\xi )\;\;.$$
we refer to the  section \ref{pow} where a numerical example shows that in 
case of bundle of perpendicular, not necessarily narrow, channels $G$ is 
differentiable and linear, i.e. that (\ref{conjj,1}) can be expected to hold 
in rectangular geometry. That can be explained by the fact that we have $n$ 
unidirectional flows with, almost, no interaction between them. In fact, 
example \ref{e,1} proves that the interaction between two flows exists, but it 
has some lower order.\\
As it was noticed in \cite{lec} the function $\Uc$ is differentiable for any 
$\xi\neq 0$. We are, therefore, in situation very close to (\ref{conjj,1}).
Moreover
\begin{equation}
\frac{\partial \Uc (\xi)_i}{\partial \xi_j}=\int_{\Yc} [w^{j}_\xi 
(y)]_i\;dy\;\;,\label{dif,1}
\end{equation}
where
\begin{eqnarray}
&&-\mu\, \mbox{div}_y \{(r-2)|e_y(w_\xi)|^{r-4}\,[e(w_\xi)\cdot e(w^j_\xi )]\, 
e_y(w_\xi)+\\
&&+|e_y(w_\xi)|^{r-2} e_y(w^{j}_\xi)\}+\nabla_y \pi^{j}_\xi  =  {\bf 
e}_j\;\;\mbox{in}\;\; \Yc\label{deri,1}\\
&&\mbox{div}_y w^j_\xi =0 \;\;\mbox{in}\;\; \Yc\nonumber\\
&&(w^j_\xi,\pi^j_\xi)\;\;\mbox{is}\;\Yc -\mbox{periodic }\nonumber\\
&&w^j_\xi =0 \;\;\mbox{on}\;S\;\;.\label{dera,1}
\end{eqnarray}
Since this problem is linear it is easy to prove the following result:
\begin{prop}
For $\xi\neq 0$ the problem (\ref{deri,1})-(\ref{dera,1}) has a unique 
solution $w^i_\xi \in V_\xi $, where 
$$V_\xi =\{\phi \in W^{1,r}_{\#}(\Yc)^n\;;\;|e(w_\xi)|^{\frac{r}{2}-1} e(\phi) 
\in L^2 (\Omega )\;\;,\;\mbox{div}\,\phi=0\;\}\;\;$$
and
$$W^{1,r}_{\#}(\Omega)=\{\phi\in W^{1,r}(\Yc)\;;\;w^j_\xi\;\mbox{is}\;\Yc 
-\mbox{periodic }\;,\;w^j_\xi 
=0 \;\mbox{on}\;S\}\;.$$
Furthermore $w^j_\xi $ satisfies the a priori estimate
\begin{equation}
|e(w^j_\xi ) |_{L^r (\Yc )} \leq C\,|e(w_\xi )|_{L^r (\Yc 
)}^{2-r}\;\;.\label{es,1}\end{equation}
\end{prop}
{\bf Proof.} We first notice that the definition of $V_\xi $ has a sense 
since, due to the (\ref{osnovna,1}), $e(w_\xi)\neq 0 $ (a.e) in $\Yc $. Now 
$V_\xi$ is the Banach space equipped by the norm
$$|\phi|_{V_\xi} = |e(\phi)|_{L^r ({\cal Y})}+ |\;|e(w_\xi)|^{\frac{r}{2}-1} 
e(\phi)|_{ L^2 
({\cal Y} )}\;\;.$$
The operator $T:V_\xi \to V_\xi'$, defined by
$$T\,\phi = -\mu\, \mbox{div}_y \{(r-2)|e_y(w_\xi)|^{r-4}\,[e(w_\xi)\cdot 
e(\phi )]\, e_y(w_\xi)+|e_y(w_\xi)|^{r-2} e_y(\phi)\} $$
is coercive since, using the H\"{o}lder's inequality, we get
\begin{eqnarray*}&\langle T\phi |\phi \rangle &\geq (r-1)\int_{\Yc} 
|e(w_\xi)|^{r-2}\,|e(\phi)|^2\\
&&\geq \frac{r-1}{2}\{|\;|e(w_\xi)|^{\frac{r}{2}-1} e(\phi)|_{ L^2 (\Omega 
)}^2+ (|e(w_\xi)|_{L^r(\Yc)})^{(r-2)} |e(\phi)|^2_{L^r (\Yc )}\}\;.
\end{eqnarray*}
Now $T\in {\cal L}(V_\xi ,V_\xi')$ is linear and coercive. Therefore it is 
bijective. The estimate (\ref{es,1}) follows from
$$\langle T\,w^j_\xi |w^j_\xi \rangle \geq 
(r-1)(|e(w_\xi)|_{L^r(\Yc)})^{(r-2)} |e(w^j_\xi)|^2_{L^r (\Yc 
)}\;\;.\;\;\clubsuit $$
\begin{lm}
The matrix $\nabla_\xi \,{\cal U} (\xi)\;,\;\xi\neq 0\;$ is symmetric.
\end{lm}
{\bf Proof.} Multiplying (\ref{deri,1}) by $w^j_\xi$ and integrating over $\Yc 
$ we obtain
\begin{eqnarray*}
&&\frac{\partial \Uc_j}{\partial \xi_i} 
(\xi)=\int_{\Yc}\{(r-2)|e_y(w_\xi)|^{r-4}\,[e(w_\xi)\cdot e(w^i_\xi )]\, 
[e_y(w_\xi)\cdot e_y(w^j_\xi)]+\\
&&+|e_y(w_\xi)|^{r-2} e_y(w^{j}_\xi)\cdot e_y (w^i_\xi )\}=\frac{\partial 
\Uc_i}{\partial \xi_j} (\xi)\;\;.\;\clubsuit\end{eqnarray*}
Now $\Uc $ is differentiable, its derivative is defined by 
(\ref{deri,1})-(\ref{dera,1}) in any point $\xi \in {\bf R}^n$ except for 
$\xi=0$  and it is symmetric. Consequently $G$ is differentiable for any 
$\xi\neq 0$ 
with
\begin{equation}
\frac{\partial G_j}{\partial \xi_i} (\xi )=(r-1)|\Uc(\xi 
)_j|^{r'-2}\frac{\partial \Uc_j}{\partial\xi_i} 
(\xi)\;\;.\label{izo,1}\end{equation}
Due to the estimate 
$$\frac{|\Uc(\xi)-\Uc (0)|}{|\xi |}\leq M\;|\xi|^{r'-2} $$
we see that $\nabla_\xi\,\Uc (0) =0 $. 
It is, therefore, not possible to linearize ${\cal U}$ in vicinity of $\xi=0$, 
as we did in the Carreau's case. Once again we have the problem with 
$\nabla_\xi\, G $ for $\xi=0 $, since (\ref{izo,1}) becomes an undetermined 
expression of the form $\frac{0}{0}$.
But, we obviously have
$$|\nabla_\xi\,G(\xi )| \leq C\;\;,$$
so that $G\in C^{0,1} ({\bf R}^n)$.\\
Unfortunately, this is as far as we go with (\ref{conjj,1}). 
The  two numerical examples found in section \ref{pow} show that, in general, 
$G$ is nonlinear. That, of course, means that it is not differentiable for 
$\xi=0 $, which can also be seen from those numerical computations.
\section{Numerical Experiments}
To illustrate the various theoretical results of this article, we performe in 
this section  numerical computations. After a short description of the 
numerical schemes used to
solve the linear and nonlinear  Stokes problem, we describe the different 
pores geometries considered. The following section (\ref{car}) is devoted to 
the  numerical experiment in case of Careau's viscosity, and finally the last 
section (\ref{pow}) is concerning the power-law viscosity.
\subsection{Numerical methods used}
We have to compute the solution of the nonlinear Stokes system 
(\ref{osnovna,1}--\ref{osnovno,1}) and for the case of the Carreau law 
viscosity, the solution of linear Stokes problems. Two difficulties have to be 
solved: the incompressibility condition and the nonlinearity in case of 
nonlinear viscosity law. An Augmented Lagrangian algorithm is used to satisfy 
the incompressibility condition $\hbox{div }u = 0$. We use a finite element 
method to discretize the Stokes problems: a $P_2^\circ$ is used for the 
velocity field, and a $P_1$ discontinuous element is used for the pressure 
field. Those elements are known to give good results with Augmented Lagrangian 
algorithm. Finnaly, the nonlinearities caused by the viscosity law are solved 
by a Newton Raphson procedure. We refer to \cite{GZ} for the description of 
the global algorithm. The porous cells are discretized by a mesh containing 
1000 triangles.
\subsection{The Porous media considered}
We consider three models of porous media to make numerical experiments: Each 
of them are characterized by the solid inclusion in the periodic unitary cell. 
As shown in figure (\ref{fig:porous model}), we consider three geometries of 
the inclusion: a cylindrical inclusion, a square inclusion, and finally, an 
ellipsoidal inclusion which gives an anisotropic porous medium. In the 
following, we refer to {\bf GEOM1} for the prous medium with circular 
inclusion, {\bf GEOM2} for this one with square inclusion and finally, {\bf 
GEOM3} for the last porous medium with ellipsoidal inclusion. 
\begin{figure}[h]  
\centerline{\includegraphics[width=5.cm]{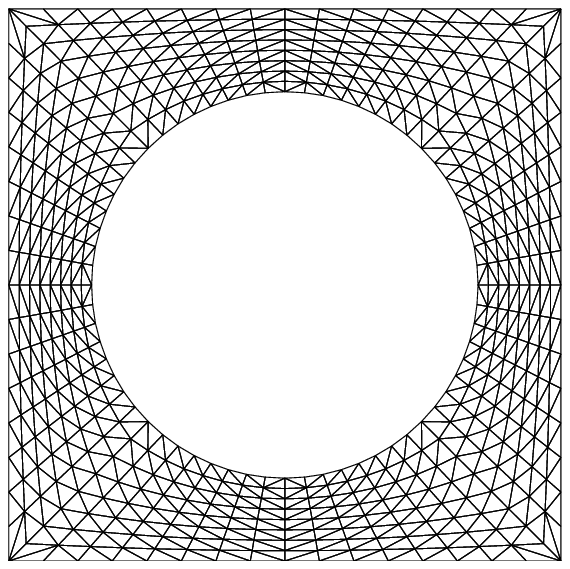} 
            \includegraphics[width=5.cm]{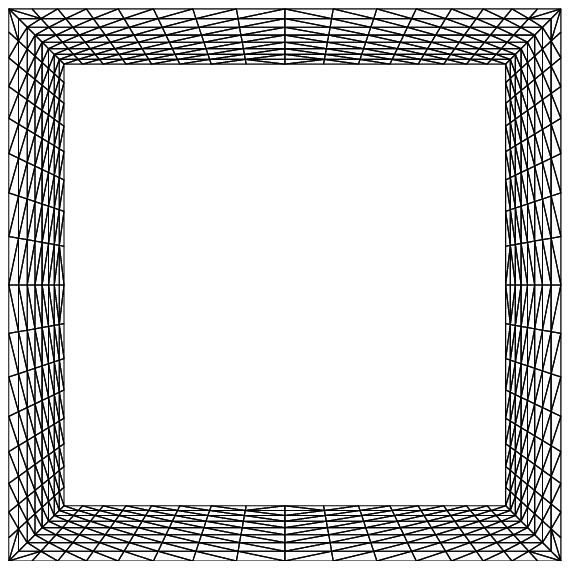} 
            \includegraphics[width=5.cm]{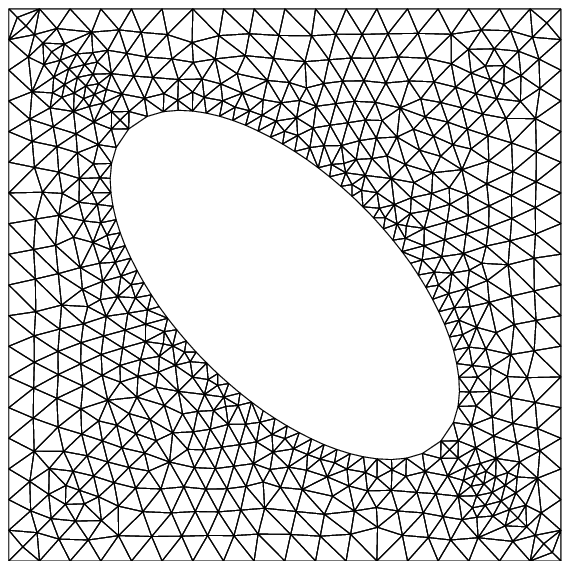}}
\caption{ Periodic unitary cell: {\bf GEOM1} on the left, {\bf GEOM2} in the 
middle and finally {\bf GEOM3} on the right of the figure.} \label{fig:porous 
model} 
\end{figure} 
\subsection{Numerical Simulations for Carreau's viscosity law}\label{car}
\begin{figure}[h] 
\centerline{\includegraphics[width=8.cm]{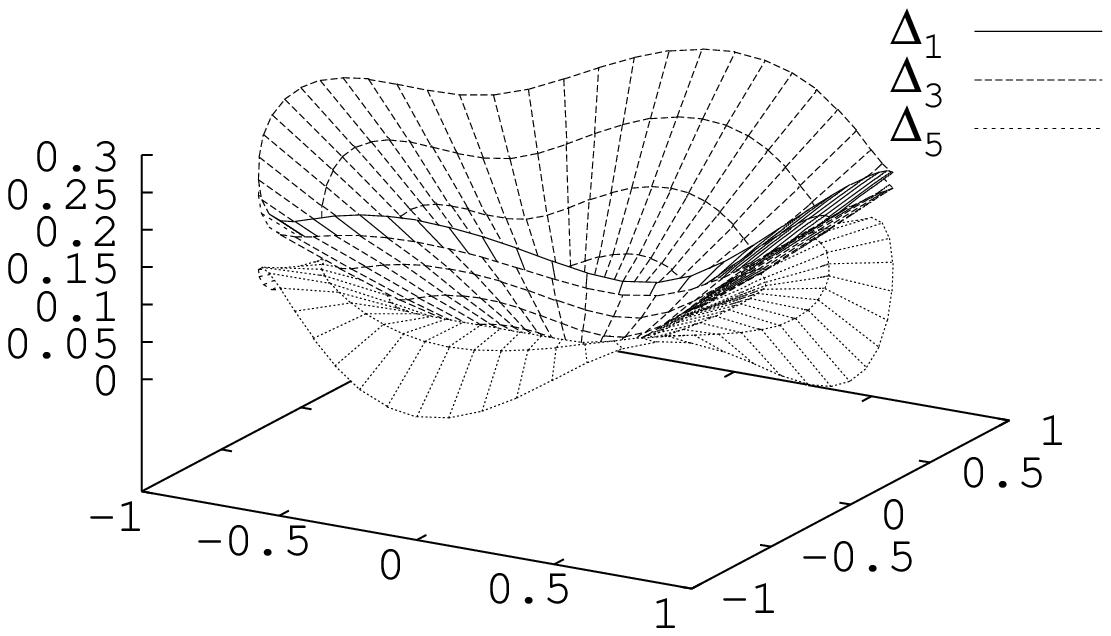}
            \includegraphics[width=8.cm]{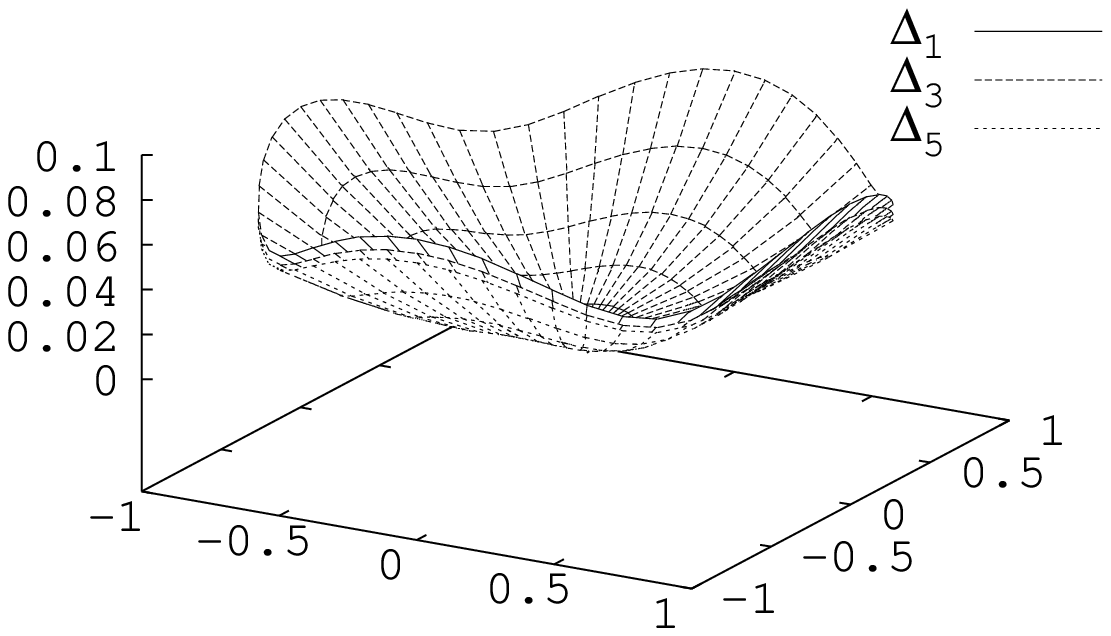}}
\caption{  Gap between ${\cal U}$ and its Taylor's expansions of order 1, 3 
and 5
 for {\bf GEOM1} on the left, and for {\bf GEOM2} on the right (for $\lambda_2 
= 100$).} \label{fig:cer-car100} 
\end{figure}

\begin{figure} 
\centerline{\includegraphics[width=8.cm]{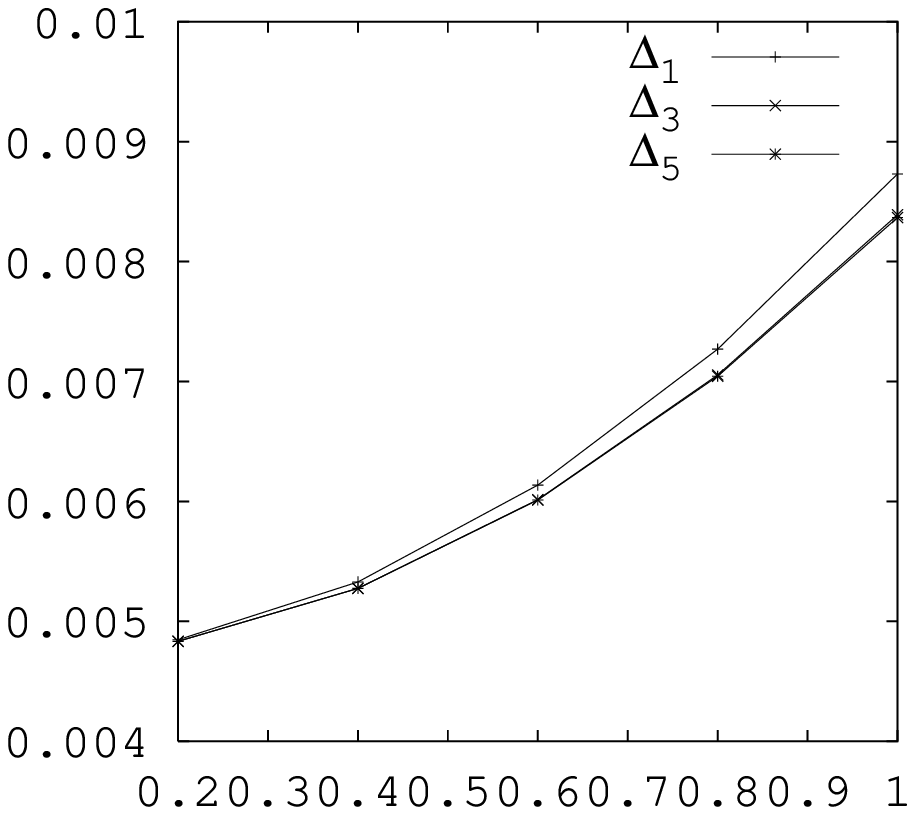}
            \includegraphics[width=8.cm]{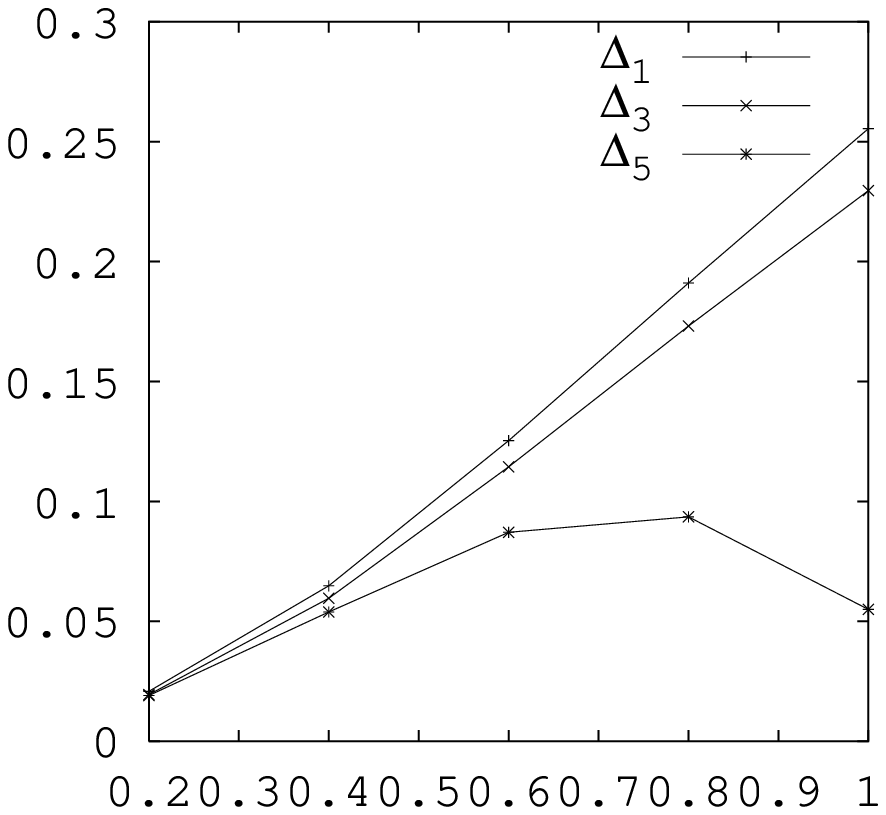}}
\caption{Gap between ${\cal U}$ and its Taylor's expansion for {\bf 
GEOM1}:$|\xi | \leq 1$, $ ({\bf e}_1 ,\xi ) = 45^\circ,\;\lambda=1$ on 
the left, $ \lambda=100 $ on the right.} 
\label{fig:cer} 
\end{figure}

\begin{figure}  
\centerline{\includegraphics[width=8.cm]{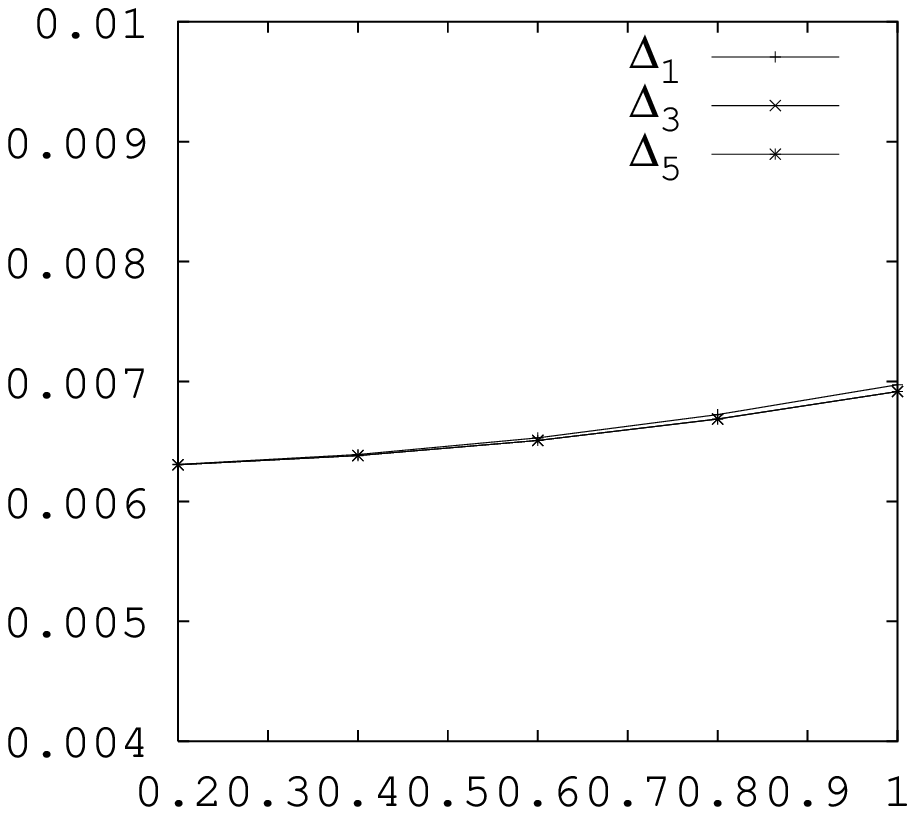}
            \includegraphics[width=8.cm]{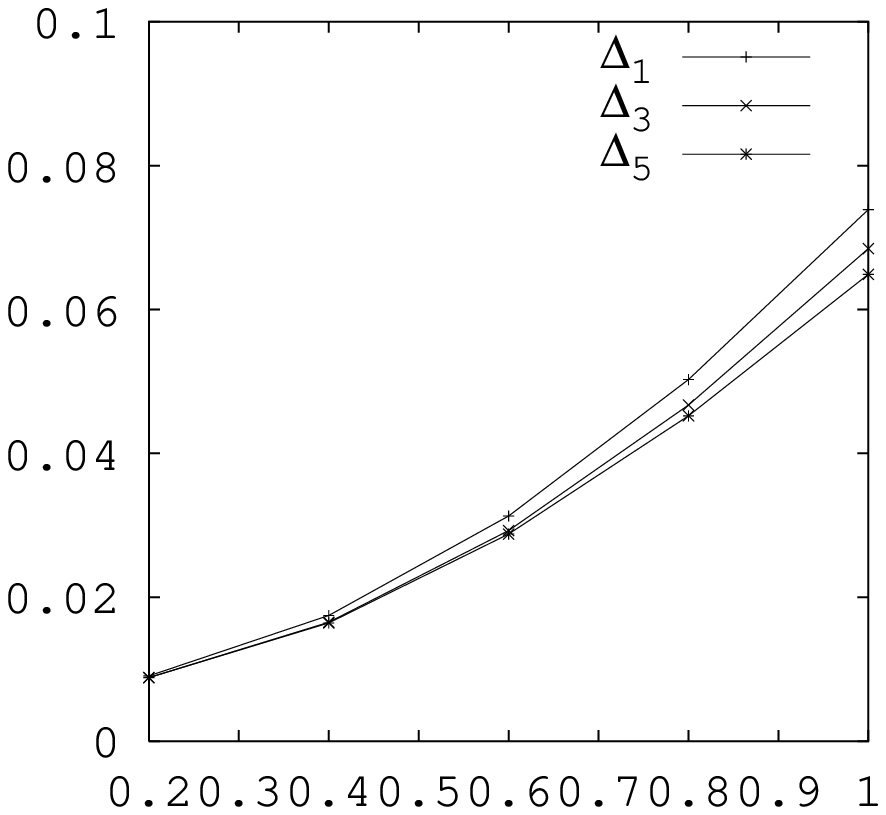}}
\caption{Gap between ${\cal U}$ and its Taylor's expansion for {\bf 
GEOM2}:$|\xi| \leq 1$, 
$ ({\bf e}_1 ,\xi ) = 45^\circ$, $\lambda=1$ on the left, 
$\lambda=100$ on the right.} \label{fig:car} 
\end{figure}

\begin{figure}[h]  
\centerline{\hspace{1.5cm}
    \includegraphics[width=7.cm]{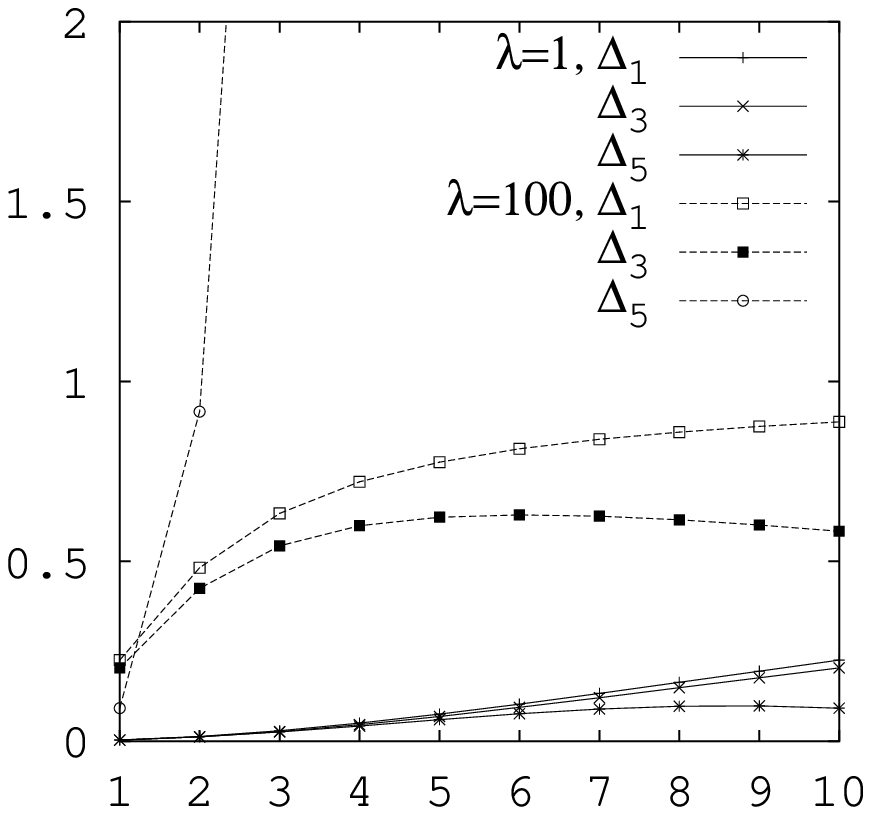}\hspace{-1.5cm}
    \includegraphics[width=7.cm]{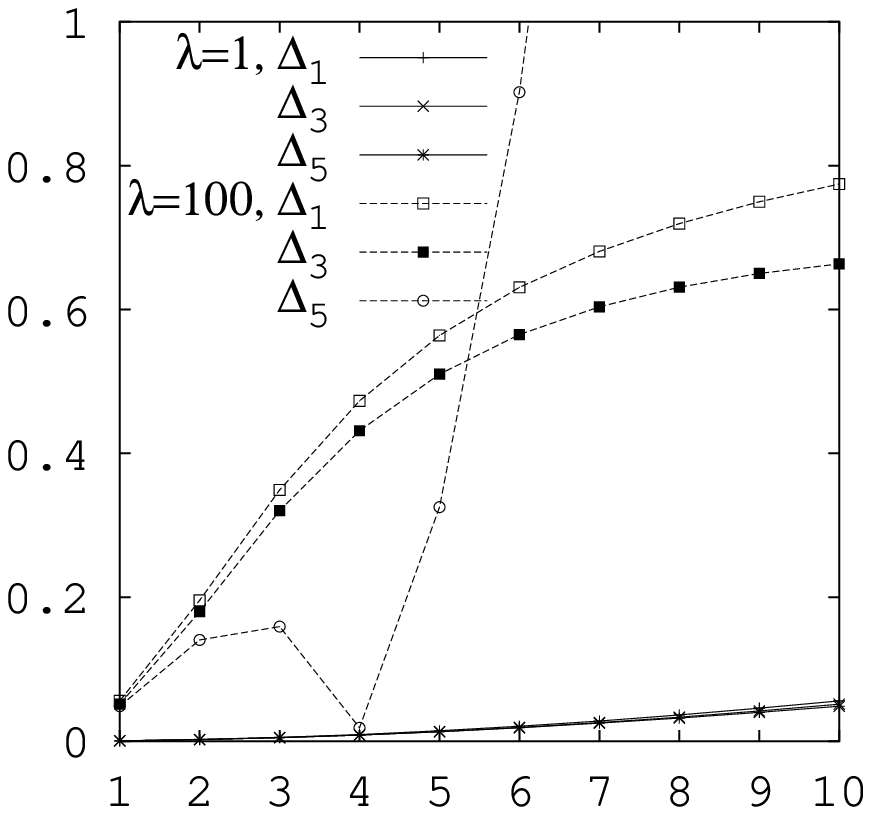}\hspace{-1.5cm}
    \includegraphics[width=7.cm]{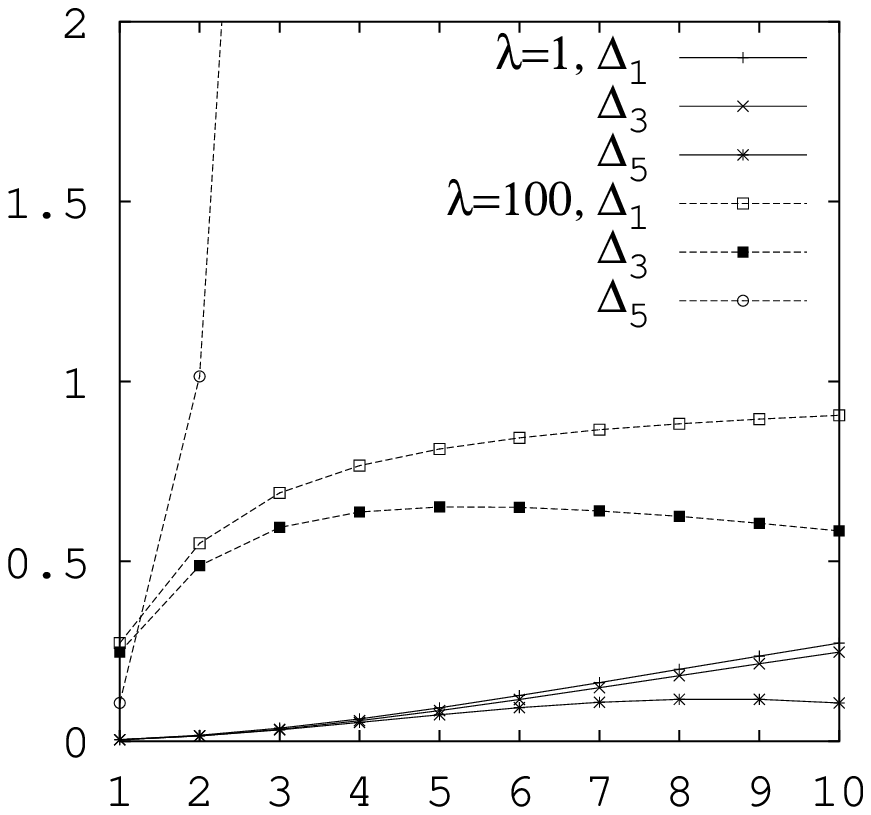}}
\caption{  Gap between ${\cal U}$ and its Taylor's development of order 1, 3 
and 5 for $|\xi| \le 10$, for {\bf GEOM1} on the left, and for {\bf GEOM2} in 
the middle and  {\bf GEOM3} on the right.} \label{fig:cer-car100_10} 
\end{figure}

In this section, we propose to solve (\ref{osnovna,1}--\ref{osnovno,1}) for a 
sampling of external force $\xi$ and its Taylor's expansion around the origin 
when viscosity is governed by a Carreau's law. At first, we describe 
parameters of the Carreau law and the numerical experiments procedure. In 
applications (see e.g. \cite{BAH}) $\eta_\infty$ is either very small compared 
to $\eta_0$ (polymer solutions) or even equal to zero (polymer melts). For 
this reason, neglect $\eta_\infty $ and consider a simplified Carreau law 
viscosity with $\eta_\infty = 0$. We chose arbitrarely $\eta_0 = 1$ and 
consider two values of $\lambda$: $\lambda_1=1$ and $\lambda_2 = 100$ to take 
into account the different comportment for a viscosity law: one near a 
constant law, and the other one far of the constant case. Finally, we consider 
a flow index $r = 1.5$:
\begin{displaymath}
\eta(e(u)) = (1+ \lambda_i |e(u)|^2)^{-.25} \quad \hbox{with } \lambda_1=1. 
\hbox{ or } \lambda_2 = 100.
\end{displaymath}
Because theoretical results are given for small $\xi$, we consider at first, 
the solution of the nonlinear problem (\ref{osnovna,1}--\ref{osnovno,1}) for 
$|\xi|\le 1$. Numerical solution are performed for a sampling $\xi_{ij}$ of 
$\xi$ defined by:
\begin{displaymath}
\xi^{ij} = {}^t(\frac{i}{n}\hbox{cos} \frac{j\pi}{4m}, \frac{i}{n}\hbox{sin} 
\frac{j\pi}{4m}), 1\le i\le n = 5, \;0\le j\le m = 8 
\end{displaymath}
extended by symmetry for $m<j<4m$ . For each value of $\xi^{ij}$ we compute 
the gap $\Delta^1_{ij}$ (resp.$\Delta^3_{ij}$, $\Delta^5_{ij}$)  between the 
solution ${\cal U}_{ij}$ of (\ref{osnovna,1}--\ref{osnovno,1}) and the 
Taylor's expansion of order $1$ (resp $3$, $5$)  where $\Delta^k_{ij}$ are 
defined by:
\begin{eqnarray}
 \Delta^1_{ij} &=& \left [ \displaystyle \sum_{k=1}^2({\cal U}_{ij} - K   
\xi^{ij})_k^2\right]^{\frac{1}{2}} \left [\displaystyle \sum_{k=1}^2 {\cal 
U}_k^2 \right ]^{-\frac{1}{2}},  \label{eq:deltaU1} \\  \Delta^3_{ij} &=& 
\left [\displaystyle \sum_{k=1}^2\left ({\cal U}_{ij} - K    \xi^{ij} - 
\frac{1}{2} \lambda (r-2)\sum_{\ell,m,n,p=1}^2  H_{\ell 
m}^{np}\xi^{ij}_\ell\xi^{ij}_m\xi^{ij}_n {\bf e}_p \right )_k^2 
\right]^{\frac{1}{2}}\left [\displaystyle \sum_{k=1}^2 {\cal U}_k^2 \right 
]^{-\frac{1}{2}} \label{eq:deltaU2}\\
  \Delta^5_{ij} &=& \left [\displaystyle \sum_{k=1}^2\left ({\cal U}_{ij} - K  
  \xi^{ij} - \frac{1}{2} \lambda (r-2)\left\{\sum_{\ell,m,n,p=1}^2  H_{\ell 
m}^{np}\xi^{ij}_\ell\xi^{ij}_m\xi^{ij}_n {\bf e}_p \right . \right. \right. 
\nonumber \\&+&\left.\left. \left. \sum_{\ell,m,n,p,q,r=1}^2  H_{\ell 
mn}^{pqr}\xi^{ij}_\ell\xi^{ij}_m\xi^{ij}_n \xi^{ij}_p\xi^{ij}_q {\bf e}_r 
\right \}\right )_k^2 \right]^{\frac{1}{2}}\left [\displaystyle \sum_{k=1}^2 
{\cal U}_k^2 \right ]^{-\frac{1}{2}} \label{eq:deltaU3}
\end{eqnarray} 
On figure (\ref{fig:cer-car100}) we plot different $\Delta_{ij}^k$ for 
$k=1,2,3$ for $\lambda_2=100$ and for the first two geometries {\bf GEOM1} and 
{\bf GEOM2}. Concerning the first geometry, one may see that we need the third 
term in the Taylor's expansion to obtain good results. But, in the case of big 
square inclusion ({\bf GEOM2}), for this range of $|\xi|$, the global 
filtration law is very close to the Darcy's law: there is no real gap between 
Darcy's law and the Taylor's expansion of any of three orders.
\noindent To underline those comportments, we have drawn on figures 
(\ref{fig:cer}-\ref{fig:car}) 
the same quantities (eq. (\ref{eq:deltaU1}--\ref{eq:deltaU3})) for a fixed 
direction of $\xi$ and 
for length $0\le |\xi | \le 1$. Let consider first figure (\ref{fig:cer}) 
concerning {\bf GEOM1}. 
In the case of $\lambda_1=1$ one may see  that for the range of considered 
$\xi$, the filtration 
law is very close to Darcy's law: the relatives errors, independently of 
order, are smaller than 
$1\%$. But in the case of $\lambda_2=100$, one needs the fifth approximation 
to obtain error 
smaller than $1\%$. Secondly, when considering figure  (\ref{fig:car}) 
concerning {\bf GEOM2}, 
one may see that independently of the value of $\lambda$, the behaviour of 
${\cal U}$ remains 
almost linear: for $\lambda_1=1$ error is smaller than $0.8\%$ and for 
$\lambda_1=100$ error is 
smaller than $8\%$. To close this  numerical study of theoretical results 
concerning Carreau's 
viscosity, we have computed solution of (\ref{osnovna,1}--\ref{osnovno,1}) and 
the different gaps 
(eq. (\ref{eq:deltaU1}--\ref{eq:deltaU3})) for $|\xi|$ larger than 1, to 
obtain numerical limit 
of validity of the Taylor's expansion, for the three geometries but only for 
one fixed direction  
of $\xi$: $\angle ({\bf e}_1,\xi ) = 60^\circ$. Those results are plotted on 
figure 
\ref{fig:cer-car100_10}. One may see that for the three geometries, Taylor's 
expansions give 
good results for $\lambda_1=1$ for the considered range of $\xi$, but in the 
case of larger 
parameter $\lambda_2=100$ one loses precision of the Taylor's expansion. More 
precisely, the 
best results are for the third order Taylor's expansion, but the fifth order 
expansion gives
 very bad results ( the rest of the Taylor expansion  of order $7$ becomes non 
negligible). 
This suggests that for large $|\xi |$ and $\lambda $ the Taylor's series 
diverge. 
\subsection{Numerical simulation for Power viscosity law}\label{pow} In this 
section, we propose to solve (\ref{osnovna,1}--\ref{osnovno,1}) for a sampling 
of $\xi$ when viscosity is governed by a Power law to illustrate the 
theoretical results of section \ref{sepl,1}. At first, we describe parameters 
of the Power law and the numerical experiments procedure. Due to the 
homogeneity property of the power law viscosity, we just consider a viscositiy 
law with  $\lambda=1$. We chose arbitrarely a flow index $r = 1.5$ :
\begin{displaymath}
\eta(e(u)) = |e(u)|^{-.5} 
\end{displaymath}
For the same reason of homogeneity, we consider only the solution of the 
nonlinear problem (\ref{osnovna,1}--\ref{osnovno,1}) for $|\xi| = 1$. 
Numerical solution are performed for a sampling $\xi^{j}$ of $\xi$ defined by:
\begin{equation}\label{sampl}
\xi^{j} = {}^t(\hbox{cos} \frac{j\pi}{4m},\hbox{sin} \frac{j\pi}{4m}),\;0\le 
j\le m = 8 
\end{equation}
and any extended by symmetry as in the Carreau's case. In two following 
sections, we will numerically show the different behaviour of the filtration 
law depending on the geometry of porous medium.
\subsubsection{Numerical results}\label{pow_1}
\begin{figure}[h]
\centerline{
       \includegraphics[width=6.5cm]{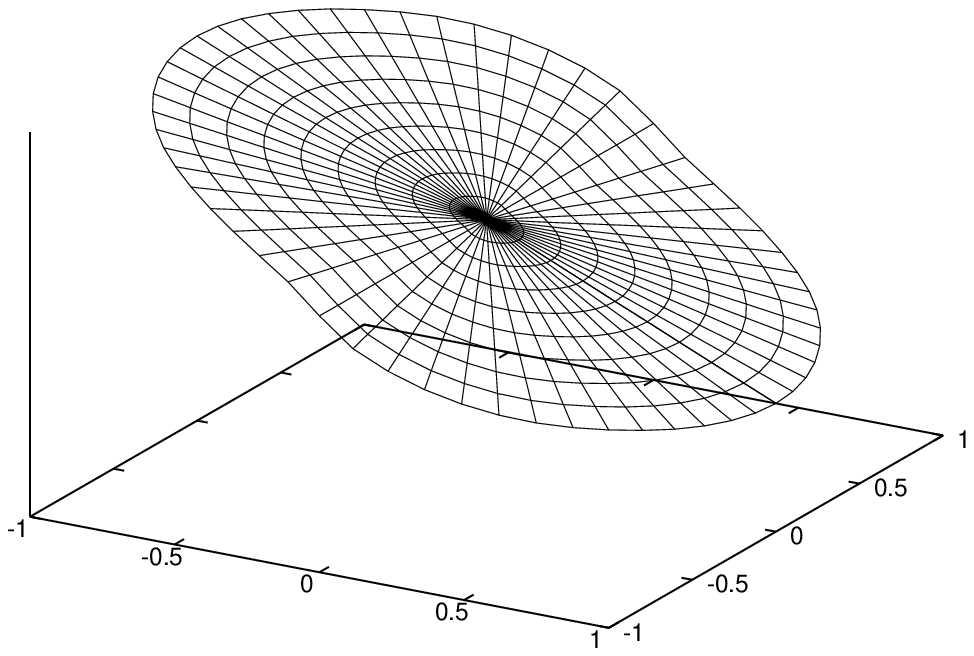}\hspace{-1.5cm}
       \includegraphics[width=6.5cm]{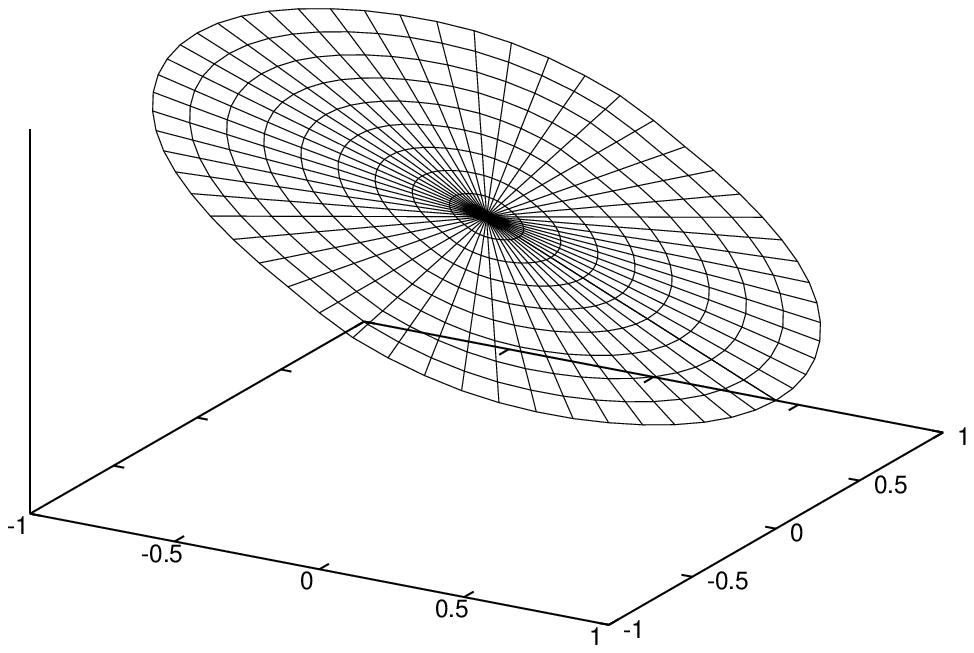}\hspace{-1.5cm}
       \includegraphics[width=6.5cm]{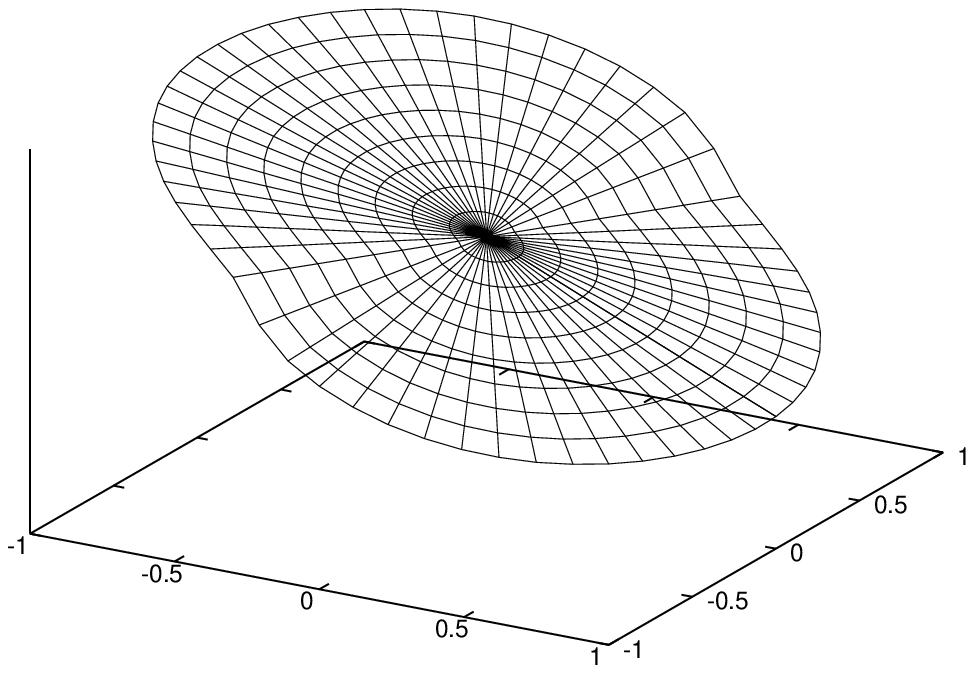}}
\caption{First component $G_1(\xi$) of the function $G$ for{ \bf GEOM1} on 
left { \bf GEOM2}
in the middle and { \bf GEOM3} on right} \label{fig:G1} 
\end{figure}

\begin{figure}[h]  
\centerline{\input{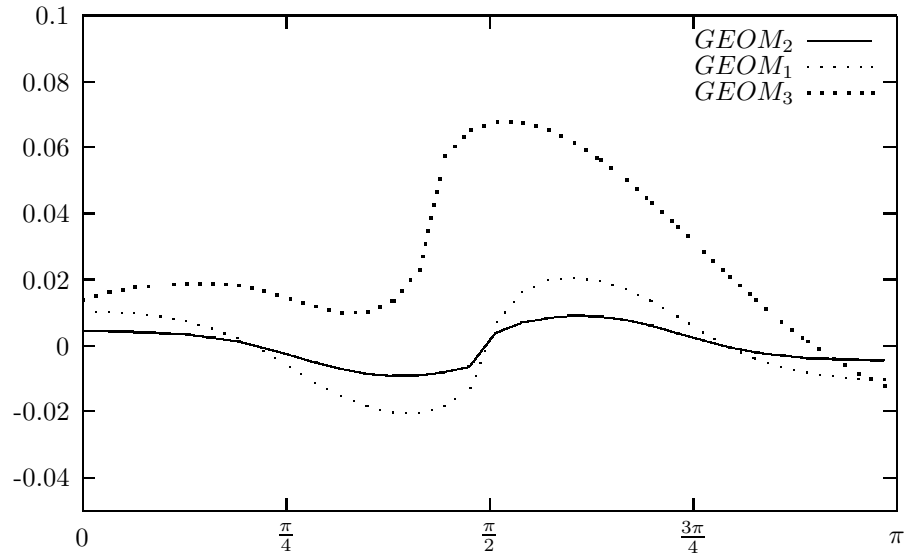}}
\caption{Pointwise gap between $G_1$ and its least square approximation 
$y=A\xi_1$ and $G_1(\xi)$ for {\bf GEOM1},{\bf GEOM2} and{\bf GEOM3}. } 
\label{fig:gap} 
\end{figure}

This section is devoted to numerical results in the case of power law to 
exhibit the linear or nonlinear comportment of the function $G(\xi)$ defined 
by its components 
$G(\xi)_i = | {\cal U}_i(\xi) |^{r-2} {\cal U}_i(\xi)$. On the figure 
(\ref{fig:G1}) we plot the first component $G_1$ of the function $G$ for the 
three geometries $\xi$ (due to the symmetry of inclusion, the comportment of 
the second component of $G$ is deduced from the first one by an appropriate 
symmetry). It may be noticed that the comportment of $G_i$ as a function of 
$\xi_i$ seems linear  for the case of the square inclusion ({\bf GEOM2}), but 
one sees  that in the case of the spherical or ellipsoidal inclusion ({ \bf 
GEOM1} and{ \bf GEOM3}) , this is not the case.
To underline this effect, we have computed the slope $A$ of the straight line 
$y = A \xi_i$ which approaches, in the least square  sense, the function $G_i$ 
for the each of the three geometries. This slope is solution of the following 
problem:
\begin{displaymath}
\inf_{B\in \re} \sum_{j=1}^N  [ G_i(\xi^j) - B \xi^j_i ]^2 \end{displaymath} 
and is given by:
\begin{displaymath}
A = \frac{\displaystyle \sum_{j=1}^N G_i(\xi^j) \xi^j_i }{ \displaystyle 
\sum_{j=1}^N(\xi^j_i)^2 }.
\end{displaymath}
Then we have computed the relative distance $\Delta$ between the straight line 
and the graph of $G_i$:
$$\Delta =\left [ 
\frac{\sum_{j=1}^N[G_i(\xi^j)-A\xi^j_i]^2}{\sum_{j=1}^N[G_i(\xi^j)]^2}\right 
]^{\frac{1}{2}}$$
In table\ref{tab1} we give the values of $A$ and $\Delta$:
\begin{table}[h]\begin{center}\footnotesize
\begin{tabular}{|c|c|c|c|}\hline
& $ GEOM_1 $ &  $ GEOM_2 $  &  $ GEOM_3 $    \\ \hline\lower.3ex\hbox{A}     & 
\lower.3ex\hbox{$-0.24083\;10^{-1} $} & \lower.3ex\hbox{$  -0.63016\;10^{-2} 
$} &  \lower.3ex\hbox{$ -0.41984\;10^{-1} $ }   \\ \hline  $\Delta$   &$ 
0.99781\;10^{-1}  $ & $ 0.44640\; 10^{-1} $ & $ 0.27976 $   \\ \hline
\end{tabular}
\end{center}
\caption{Slope $A$ of the straight line and relative distance $\Delta$ between 
this line and $G_1$ for the three geometries}
\label{tab1}
\end{table}
\noindent One sees that $\Delta $ is the smallest for the square inclusion  
({\bf GEOM2}), (of order $4 \%$), and is no negligible for the other 
geometries ($10\%$ for the spherical inclusion and $27\%$ for the ellipsoidal 
inclusion). Finally, we have drawn on figure (\ref{fig:gap}), the pointwise 
distance between the straight line $y=A\xi_1$ and the function $G(\xi)_1 = | 
{\cal U}_1(\xi) |^{r-2} {\cal U}_1(\xi)$ for the three geometries. One may see 
that the gap between the straight line and the function $G_1$ is maximal at 
the point where ${\cal U}_1$ is changing  sign (for angle $\pi/2$ for  { \bf 
GEOM1} and {\bf GEOM2}). This effect may be understood as the 
non-differentiability of ${\cal U}$ in the origin.

\vspace{.5cm}

\noindent{\bf Acknowledgement.} This paper was written during the stay of 
Eduard Maru\v{s}i\'{c}-Paloka at Universit\'{e} de Saint-Etienne during the 
June of 99.


\end{document}